\documentclass[aop,seceqn,MSNbibl,citesort,dvips]{arximspdf}

\doi{10.1214/09-AOP513}
\volume{38}
\issue{4}
\pubyear{2010}
\firstpage{1401}
\lastpage{1443}

\makeatletter

\newproclaim{definition}{Definition}[section]
\newproclaim{remark}{Remark}
\newtheorem{lemma}{Lemma}
\newtheorem{proposition}{Proposition}
\newtheorem{theorem}{Theorem}

\newcommand{\dd}{{d}}
\newcommand{\ii}{{i}}
\newcommand{\nnorm}{|\!|\!|}

\newproclaim{Step}{Step}

\makeatother

\begin{document}
\begin{frontmatter}

\title{On ergodicity of some Markov processes}
\runtitle{Ergodicity of invariant measures}

\begin{aug}
\author[A]{\fnms{Tomasz} \snm{Komorowski}\corref{}\thanksref{t1,t4}\ead[label=e1]{komorow@hektor.umcs.edu.pl}},
\author[B]{\fnms{Szymon} \snm{Peszat}\thanksref{t2,t4}\ead[label=e2]{napeszat@cyf-kr.edu.pl}} and
\author[C]{\fnms{Tomasz} \snm{Szarek}\thanksref{t3,t4}\ead[label=e3]{szarek@intertele.pl}}
\runauthor{T. Komorowski, S. Peszat and T. Szarek}
\affiliation{Maria Curie-Sk\l odowska University and Polish Academy of
Sciences,
Polish~Academy of Sciences and University of Gda\'nsk}
\dedicated{To the memory of Andrzej Lasota (1932--2006)}
\address[A]{T. Komorowski\\
Institute of Mathematics\\
Maria Curie-Sk\l odowska University\\
Pl. M. Curie-Sk\l odowskiej 1, 20-031 Lublin\\
Poland\\
and\\
Institute of Mathematics\\
Polish Academy of Sciences\\
\'Sniadeckich, 8, 00-956, Warsaw\\
Poland\\
\printead{e1}}
\address[B]{S. Peszat\\
Institute of Mathematics\\
Polish Academy of Sciences\\
\'Sw. Tomasza 30/7, 31-027 Krak\'ow\\
Poland\\
\printead{e2}}
\address[C]{T. Szarek\\
Institute of Mathematics\\
University of Gda\'nsk\\
Wita Stwosza 57, 80-952 Gda\'nsk\\
Poland\\
\printead{e3}}
\end{aug}

\thankstext{t1}{Supported in part by Polish Ministry of Science and
Higher Education Grant N20104531.}

\thankstext{t2}{Supported in part by Polish Ministry of Science and
Higher Education Grant PO3A03429.}

\thankstext{t3}{Supported in part by Polish Ministry of Science and
Higher Education Grant N201 0211 33.}

\thankstext{t4}{Supported by EC FP6 Marie Curie ToK programme SPADE2,
MTKD-CT-2004-014508 and
Polish MNiSW SPB-M.}

\received{\smonth{2} \syear{2009}}
\revised{\smonth{10} \syear{2009}}

%
\begin{abstract}
We formulate a criterion for the existence and uniqueness of an
invariant measure for a Markov process taking values in a Polish phase
space. In addition, weak-$^*$ ergodicity, that is, the weak convergence
of the ergodic averages of the laws of the process starting from any
initial distribution, is established. The principal assumptions are the
existence of a lower bound for the ergodic averages of the transition
probability function and its local uniform continuity. The latter is
called the \textit{e-property}. The general result is applied to
solutions of some stochastic evolution equations in Hilbert spaces. As
an example, we consider an evolution equation whose solution
describes the Lagrangian observations of the velocity field in the
passive tracer model. The weak-$^*$ mean ergodicity of the
corresponding invariant measure is used to derive the law of large
numbers for the trajectory of a tracer.
\end{abstract}

%

\end{frontmatter}

\section{Introduction}

The lower bound technique is a useful tool in the ergodic
theory of Markov processes. It has been used by Doeblin (see
\cite{doeblin}) to show mixing of a Markov chain whose transition
probabilities possess a uniform lower bound. A somewhat different
approach, relying on the analysis of the operator dual
to the transition probability, has been applied by Lasota and Yorke,
see, for instance, \mbox{\cite{lasotamackey,LY}}. For example, in
\cite{LY}, they show that the existence of a lower bound for the
iterates of the Frobenius--Perron operator (that corresponds to a
piecewise monotonic transformation on the unit interval) implies
the existence of a stationary distribution for the deterministic
Markov chain describing the iterates of the transformation. In
fact, the invariant measure is then unique in the class of
measures that are absolutely continuous with respect to one-dimensional
Lebesgue measure. It is also statistically stable, that is, the
law of the chain, starting from any initial distribution that is
absolutely continuous, converges to the invariant measure in the
total variation metric. This technique has been extended to more
general Markov chains, including those which correspond to iterated
function systems; see, for example, \cite{LYy}. However, most of the
existing results are formulated for Markov chains taking values in
finite-dimensional spaces; see, for example, \cite{zaharpol} for a
review of the topic.

Generally speaking, the lower bound technique which we have in mind
involves deriving ergodic properties of the Markov process from the
fact that there exists a ``small'' set in
the state space. For instance, it could be compact, such that the time
averages of the mass of the process are concentrated over that set
for all sufficiently large times. If this set is compact,
then one can deduce the existence of an invariant probability measure
without much difficulty.

The question of extending the lower bound
technique to Markov processes taking values in Polish spaces that
are not locally compact is quite a delicate matter. This situation
typically occurs for processes that are solutions of stochastic
partial differential equations (SPDEs). The value of the process is
then usually an
element of an
infinite-dimensional Hilbert or Banach space. We stress here that
to prove the existence of a stationary measure, it is not enough
only to ensure the lower bound on the transition probability
over some ``thin'' set. One can show (see the counterexample
provided in \cite{Sz}) that even if the mass of the process
contained in any neighborhood of a given point is separated from
zero for all times, an invariant measure may fail to exist. In fact,
some general results concerning the existence of an invariant
measure and its statistical stability for a discrete-time Markov
chain have been formulated in \cite{Sz}; see Theorems 3.1--3.3.

In
the present paper, we are concerned with the question of finding a
criterion for the existence of a unique, invariant, ergodic probability
measure for a continuous-time Feller Markov process
$(Z(t))_{t\ge0}$ taking values in a Polish space $\mathcal X$; see Theorems
\ref{Mtheorem} and \ref{Theorem} below. Suppose that $(P_t)_{t\ge0}$
is its
transition probability semigroup. In our first
result (see Theorem \ref{Mtheorem}), we show that there exists a
unique, invariant probability measure for the process, provided
that for any Lipschitz, bounded function~$\psi$, the family of
functions $(P_t\psi)_{t\ge0}$ is uniformly continuous at any point
of $\mathcal X$ (we call this the \textit{e-property} of the semigroup)
and there
exists $z\in\mathcal X$ such that for any $\delta>0$,
%
%
\begin{equation}
\label{intro-1}
\liminf_{T\to+\infty}\frac{1}{T}
\int_0^TP_t\mathbf{1}_{B(z,\delta)}(x)\,\dd t>0 \qquad \forall x\in
\mathcal X.
\end{equation}
Here, $B(z,\delta)$ denotes the ball in $\mathcal X$ centered at $z$ with
radius $\delta$.
Observe that, in contrast to the Doeblin condition, we do not
require that the lower bound in (\ref{intro-1}) is uniform in the
state variable $x$.
If some conditions on uniformity over bounded sets are added [see
(\ref{Th1}) and
(\ref{Th2}) below], then one can also deduce the
stability of the ergodic averages corresponding to
$(Z(t))_{t\ge0}$; see Theorem \ref{Theorem}.
We call this, after \cite{zaharpol}, \textit{weak-$^*$ mean
ergodicity.}

This general result is applied to solutions of stochastic
evolution equations in Hilbert spaces. In Theorem
\ref{TGeneral}, we show the uniqueness and ergodicity of an invariant
measure, provided that the transition semigroup has the e-property and the
(deterministic) semi-dynamical system corresponding to the equation
without the noise has an attractor which admits a unique invariant measure.
This is a natural generalization of the results known for
so-called \textit{dissipative systems}; see, for example, \cite{DZ2}.

A different approach to proving the uniqueness of an invariant measure for
a stochastic evolution equation is based on the strong Feller property
of the transition semigroup (see \cite{DZ2,Ha,EH} and \cite{PZ1}) or, in a more refined
form, on \textit{the asymptotic strong Feller property} (see
\cite{HM,HM1,M}). In our Theorem \ref{TGeneral}, we
do not
require either of these properties of the
corresponding semigroup. Roughly speaking, we assume: (1) the existence
of a global compact attractor for the system without the noise
[hypothesis (i)];
(2) the existence of a Lyapunov function [hypothesis~(ii)]; (3) some
form of stochastic stability
of the system after the noise is added [hypothesis~(iii)];
(4) the e-property (see Section \ref{sec2}).
This allows us to
show lower bounds for the transition probabilities and
then use Theorems \ref{Mtheorem} and \ref{Theorem}.

As an application of Theorem \ref{TGeneral}, we consider, in Sections
\ref{Observable-section} and \ref{sec5}, the Lagrangian
observation process corresponding to the passive tracer model
$\dot{x}(t)=V(t,x(t))$, where $V(t,x)$ is a time--space stationary
random, Gaussian and Markovian velocity field. One can show that
when the field is sufficiently regular [see (\ref{H1})], the process
${\mathcal{Z}}(t):=V(t,x(t)+\cdot)$ is a solution of a certain
evolution equation in a Hilbert space; see (\ref{Zet}) below. With
the help of the technique developed by Hairer and Mattingly
\cite{HM} (see also \cite{EM} and \cite{ks}), we verify the
assumptions of Theorem \ref{TGeneral} when $V(t,x)$ is periodic in
the $x$ variable and satisfies a mixing hypothesis in the temporal
variable; see (\ref{H2}). The latter reflects, physically, quite a
natural assumption that the mixing time for the velocity field gets
shorter on smaller spatial scales. As a consequence of the
statistical stability property of the ergodic invariant measure for
the Lagrangian velocity $({\mathcal{Z}}(t))_{t\ge0}$, we obtain the
weak law of large numbers for the passive tracer model in a
compressible environment; see Theorem \ref{TTracer}. It generalizes
the corresponding result that holds in the incompressible case,
which can be easily deduced due to the fact that the invariant
measure is known explicitly in that situation; see
\cite{portstone}.

\section{Main results}
\label{sec2}

Let $(\mathcal X,\rho)$ be a Polish metric space.
Let $\mathcal B(\mathcal X)$ be the space of all Borel subsets of
$\mathcal X$ and let
$B_b(\mathcal X)$ [resp., $C_b(\mathcal X)$] be the Banach space of all bounded,
measurable (resp., continuous) functions
on $\mathcal X$ equipped with the supremum norm $\|\cdot\|_\infty$.
We denote by $\operatorname{Lip}_b (\mathcal X)$ the space of all bounded
Lipschitz continuous functions on $\mathcal X$. Denote by
\[
\operatorname{Lip} (f):= \sup_{x\not=y} \frac{|f(x)-f(y)|}{\rho(x,y)}
\]
the smallest Lipschitz constant of $f$.

Let $(P_t)_{t\ge0}$ be the transition semigroup of a Markov family
$Z=((Z^x(t))_{ t\ge0}$, $x\in\mathcal X)$ taking values in $\mathcal X$.
Throughout this paper, we shall assume that the
semigroup $(P_t)_{t\ge0}$ is \textit{Feller}, that is,
$P_t(C_b(\mathcal X
))\subset C_b(\mathcal X)$. We shall also assume that the Markov family is
stochastically
continuous, which implies that $\lim_{t\to0+}P_t\psi(x)=\psi(x)$
for all $x\in
\mathcal X$ and $\psi\in C_b(\mathcal X)$.
%
%
\begin{definition}
\label{def3.1} We say that a transition semigroup $(P_t)_{t\ge0}$ has
\textit{the e}-\textit{prop\-erty} if the family of functions $(P_t\psi)_{t\ge0}$ is
equicontinuous at every point $x$ of $\mathcal X$ for any bounded and
Lipschitz continuous function $\psi$, that is, if
\[
\forall\psi\in\operatorname{Lip}_b
(\mathcal X), x\in\mathcal X, \varepsilon>0 ,\exists\delta>0
\]
such
that
\[
\forall z\in B(x,\delta), t\ge0,\qquad |P_t\psi(x)-P_t\psi
(z)|<\varepsilon.
\]
$(Z^x(t))_{ t\ge0}$ is then called an \textit{e-process}.
\end{definition}

An e-process is an extension to
continuous time of the notion of an e-chain introduced
in Section 6.4 of \cite{MT}.

Given
$B\in\mathcal B(\mathcal X)$, we denote by ${\mathcal M}_1(B)$ the
space of all
probability Borel measures on $B$. For brevity, we write $\mathcal
M_1$ instead of ${\mathcal M}_1(\mathcal X)$. Let $(P_t^*)_{t\ge0}$ be
the dual
semigroup defined on $\mathcal M_1$ by the formula $P_t^*\mu
(B):=\int_{\mathcal X} P_t\mathbf{1}_B \,\dd\mu$ for $ B\in\mathcal
B(\mathcal X)$. Recall that $\mu_*\in\mathcal M_1$ is \textit
{invariant} for
the semigroup $(P_t)_{t\ge0}$ [or the Markov family $(Z^x(t))_{ t\ge
0}$] if $P_t^*\mu_*=\mu_*$
for all $t\ge0$.

For a given $T>0$ and $\mu\in{\mathcal M}_1$, define
$Q^T\mu:=T^{-1}\int_0^T P_s^*\mu\,\dd s$. We write $Q^T(x,\cdot)$ in
the particular case when $\mu=\delta_x$. Let
%
%
\begin{equation}\label{Pr1}
\mathcal T:= \{x\in\mathcal X\dvtx\mbox{the family of measures
$ (Q^T(x) )_{T\ge0}$ is tight} \}.
\end{equation}

Denote by $B(z,\delta)$ the ball in $\mathcal X$ with center at $z$ and
radius $\delta$, and by ``$\mathrm{w}$-$\lim$'' the limit in
the sense of weak convergence of measures. The proof of the
following result is given in Section \ref{sec2.3}.
\begin{theorem}\label{Mtheorem}
Assume that $(P_t)_{t\ge0}$ has the e-property and that there exists
$z\in
\mathcal X$ such that for every $\delta>0$ and $x\in\mathcal X$,
%
%
\begin{equation}\label{Mt1}
\liminf_{T\uparrow\infty}Q^T(x,B(z,\delta))>0.
\end{equation}
The semigroup then admits a unique, invariant probability measure
$\mu_*$. Moreover,
%
%
\begin{equation}\label{Mt2}
\mathop{\mathrm{w}\mbox{-}\mathrm{lim}}_{T\uparrow\infty}Q^T\nu= \mu_*
\end{equation}
for any $\nu\in{\mathcal M}_1$ that is supported in ${\mathcal T}$.
\end{theorem}
\begin{remark}
We remark here that the set ${\mathcal T}$ may not be the entire space
$\mathcal X$. This issue is investigated more closely in \cite{SSU}. Among
other results, it is shown there that if the semigroup $(P_t)_{t\ge0}$
satisfies the assumptions of Theorem \ref{Mtheorem}, then the set
${\mathcal T}$ is closed. Below, we present an elementary example of a
semigroup satisfying the assumptions of the above theorem, for which
${\mathcal T}\not=\mathcal X$.
Let $\mathcal X=(-\infty, -1]\cup[1, +\infty)$, $T(x):=-(x+1)/2-1$ for
$x\in\mathcal X$ and let $P\dvtx \mathcal X\times\mathcal B(\mathcal
X)\to[0, 1]$ be the
transition function defined by the formula
\[
P (x, \cdot)=
\cases{
\bigl(1-\exp(-1/x^{2})\bigr)\delta_{-x}(\cdot) +\exp
(-1/x^{2})\delta_{x+1}(\cdot),&\quad for $x\ge1$,\cr
\delta_{T(x)}(\cdot),&\quad for $x\le-1$.}
\]
Define the Markov operator $P\dvtx B_b(\mathcal X)\to B_b(\mathcal X)$
corresponding to
$P(\cdot,\cdot)$, that is,
\[
Pf(x)=\int_{\mathcal X} f(y)P (x, \dd y) \qquad\mbox{for $f\in
B_b(\mathcal X)$}.
\]
Finally, let $(P_t)_{t\ge0}$ be the semigroup given by the formula
%
%
\begin{equation}
\label{P-t-e}
P_tf=\sum_{n=0}^{\infty} e^{-t}\frac{t^n}{n!} P^nf \qquad\mbox
{for $t\ge0$}.
\end{equation}
It is obvious that the semigroup is Feller.

We check that $(P_t)_{t\ge0}$ satisfies the assumptions of Theorem
\ref{Mtheorem} and that $\mathcal T=(-\infty, -1]$. Let $z:=-1$.
Since, for every $x\in\mathcal X$ and $\delta>0$,
\[
\liminf_{t\to+\infty} P_t^*\delta_x (B(z, \delta))\ge1-\exp(-1/x^{2}),
\]
condition (\ref{Mt1}) is satisfied.

To prove the e-property, it is enough to show that for any
$f\in\operatorname{Lip}_b(\mathcal X)$,
%
%
\begin{equation}\label{unfcont}
{\lim_{y\to x}\sup_{n\ge1}} |P^n f(x)-P^n f(y)|=0 \qquad\forall
x\in\mathcal X.
\end{equation}
If $x\le-1$, then condition (\ref{unfcont}) obviously holds. We may
therefore assume that $x\ge1$.
Observe that
\[
P^n f(x)=
\sum_{k=0}^{n-1} f\bigl(T^{n-1-k}(-x-k)\bigr)G_k(x)+H_n(x)f(x+n),\qquad n\ge1,
\]
where $H_n(x):=\prod_{j=0}^{n-1}\exp(-(x+j)^{-2})$ and
$G_k(x):=[1-\exp(-(x+k)^{-2})]\times\break H_k(x)$.
Here, we interpret $\prod_{j=0}^{-1}$ as equal to $1$.
After straightforward calculations, we obtain that for $1\le x\le y$,
we have
\begin{eqnarray*}
&&
|P^n f(x)-P^n f(y)|\\
&&\qquad\le\operatorname{Lip}(f)(y-x)+\|f\|_{\infty}
\Biggl(\sum_{k=0}^{n-2} \int_x^y |G_k'(\xi)|\,\dd\xi+\int_x^y
|H_n'(\xi)|\,\dd\xi\Biggr).
\end{eqnarray*}
%
Condition (\ref{unfcont}) follows from the fact that $\sum
_{k=0}^{n-2} |G_k'(\xi)|$ and $H_n'(\xi)$ are uniformly convergent on
$[1,+\infty)$.

Finally, it can be seen from (\ref{P-t-e}) that for any $R>0$ and
$x\ge1$, we have
\[
\liminf_{t\to+\infty} P_t^*\delta_x(B^c(0,R))\ge\lim_{n\to
+\infty}H_n(x)>0,
\]
which proves that $x\notin\mathcal T$.
\end{remark}

Following \cite{zaharpol}, page 95, we introduce the notion of
weak-$^*$ mean ergodicity.
\begin{definition}
A semigroup $(P_t)_{t\ge0}$ is called \textit{weak-$^*$ mean ergodic}
if there
exists a measure $\mu_*\in{\mathcal M}_1$ such that
%
%
\begin{equation} \label{add-1}
\mathop{\mathrm{w}\mbox{-}\mathrm{lim}}_{T\uparrow\infty} Q^T\nu
=\mu_* \qquad\forall\nu\in{\mathcal M}_1.
\end{equation}
\end{definition}
\begin{remark}
In some important cases, it is easy to show
that $\mathcal T=\mathcal X$. For example, if $(Z^x(t))_{t\ge0}$ is
given by
a stochastic
evolution equation in a Hilbert space $\mathcal X$, then it is enough to
show that there exist a compactly embedded space $\mathcal
V\hookrightarrow
\mathcal X$ and a locally bounded, measurable function $\Phi
\dvtx[0,+\infty)\to
[0,+\infty)$ that satisfies
$\lim_{R\to+\infty}\Phi(R)=+\infty$ such that
\[
\forall x\in\mathcal X\ \exists T_0\ge0\qquad \sup_{t\ge T_0}
\mathbb{E}
\Phi(\|Z^x(t)\|_{\mathcal V})<\infty.
\]
Clearly, if $\mathcal T=\mathcal X$, then the assumptions of Theorem
\ref{Mtheorem} guarantee weak-$^*$ mean ergodicity. In Theorem
\ref{Theorem} below, the weak-$^*$ mean ergodicity is deduced from
a version of (\ref{Mt1}) that holds uniformly on bounded sets.
\end{remark}
\begin{remark}
Of course, (\ref{add-1}) implies uniqueness of
invariant measure for $(P_t)_{t\ge0}$. Moreover, for any stochastically
continuous Feller semigroup $(P_t)_{t\ge0}$, its weak-$^*$ mean ergodicity
also implies ergodicity of $\mu_*$, that is, that any
Borel set $B$ which satisfies $P_t\mathbf{1}_B=\mathbf{1}_B$,
$\mu_*$-a.s. for all $t\ge0$, must be $\mu_*$-trivial. This can be
seen from, for instance, part (iv) of Theorem 3.2.4 of \cite{DZ2}.
\end{remark}
\begin{remark}
Note that condition (\ref{add-1}) is equivalent to \textit{every
point of $\mathcal X$ being generic}, in the sense of \cite{furstenberg},
that is,
%
%
\begin{equation}\label{add-2}
\mathop{\mathrm{w}\mbox{-}\mathrm{lim}}_{T\uparrow\infty}  Q^T(x,\cdot)
=\mu_*\qquad \forall x\in\mathcal X.
\end{equation}
Indeed, (\ref{add-1}) obviously implies (\ref{add-2}) since it
suffices to take $\nu=\delta_x$, $x\in\mathcal X$. Conversely, assuming
(\ref{add-2}), we can write, for any $\nu\in{\mathcal M}_1$ and
$\psi\in C_b(\mathcal X)$,
\[
\lim_{T\uparrow\infty}\int_\mathcal X
\psi(x)Q_T\nu(\dd
x)=\lim_{T\uparrow\infty}\int_\mathcal X \frac{1}{T}\int_0^T
P_s \psi(x)\,\dd s\, \nu(\dd x)\stackrel{(\fontsize{8.36}{10}\selectfont{\mbox{\ref{add-2}}})}{=}\int_\mathcal X
\psi(x)\mu_*(\dd x)
\]
and (\ref{add-1}) follows.
\end{remark}

%
%
%

The proof of the following result is given in Section \ref{sec2.4}.
\begin{theorem}\label{Theorem}
Let $({P}_t)_{t\ge0}$ satisfy the assumptions of Theorem \ref{Mtheorem}.
Assume, also, that there exists $z\in X$ such that for every bounded
set $A$ and $\delta>0$, we have 
%
%
\begin{equation}\label{Th1}
\inf_{x\in A}
\liminf_{T\to+\infty}Q^T(x,B(z,\delta))>0.
\end{equation}
Suppose, further, that for every $\varepsilon>0$ and $x\in X$, there
exists a bounded Borel set $D\subset X$ such that
%
%
\begin{equation}\label{Th2}
\liminf_{T\to+\infty}Q^T(x,D)>1-\varepsilon.
\end{equation}
Then,
besides the existence of a unique invariant measure $\mu_*$ for
$(P_t)_{t\ge0}$,
the following are true:

(1) the semigroup $(P_t)_{t\ge0}$ is weak-$^*$ mean ergodic;

(2) for any $\psi\in\operatorname{Lip}_b(\mathcal X)$ and $\mu\in
{\mathcal M
}_1$, the weak law of large numbers
holds, that is,
%
%
\begin{equation}
\label{090601}
\mathop{{\mathbb{P}}_{\mu} \mbox{-}\mathrm{lim}}_{T\to+\infty
}\frac
{1}{T}\int_0^T\psi(Z(t))\,\dd
t=\int_{\mathcal X} \psi\,\dd\mu_*.
\end{equation}
Here, $(Z(t))_{t\ge0}$ is the Markov process that corresponds to the
given semigroup, whose initial distribution is $\mu$ and whose path measure
is ${\mathbb{P}}_\mu$. The convergence
takes place in ${\mathbb{P}}_\mu$ probability.
\end{theorem}

Using Theorems \ref{Mtheorem} and \ref{Theorem}, we establish the
weak-$^*$ mean
ergodicity for the family defined by the stochastic evolution
equation
%
%
\begin{equation}\label{E1}
\dd Z(t)= \bigl(AZ(t) +F(Z(t)) \bigr)\,\dd t+ R\,\dd W(t).
\end{equation}
Here, $\mathcal X$ is a real, separable Hilbert space, $A$ is
the generator of a $C_0$-semigroup $S=(S(t))_{t\ge0}$ acting on
$\mathcal X
$, $F$
maps (not necessarily continuously) $D(F)\subset\mathcal X$ into
$\mathcal X$, $R$ is a bounded linear operator from another Hilbert space
$\mathcal H$ to $\mathcal X$ and $W=(W(t))_{t\ge0}$ is a cylindrical Wiener
process on $\mathcal H$ defined over a certain filtered probability space
$(\Omega,{\mathcal F}, ({\mathcal F}_t)_{t\ge0},\mathbb{P})$.

Let $Z_0$ be an $\mathcal F_0$-measurable random variable. By a
solution of
(\ref{E1}) starting from $Z_0$, we mean a solution to the stochastic
integral equation (the so-called \textit{mild solution})
\[
Z(t)= S(t)Z_0 +\int_0^t S(t-s)F(Z(s))\,\dd s + \int_0^t S(t-s)R\,\dd
W(s),\qquad t\ge0
\]
(see, e.g., \cite{DZ1}), where the stochastic integral
appearing on the right-hand side is understood in the sense of
It\^o. We suppose that for every $x\in\mathcal X$, there is a unique mild
solution $Z^x=(Z^x_t)_{ t\ge0}$ of (\ref{E1}) starting from $x$
and that (\ref{E1}) defines a Markov family in that way. We assume
that for any $x\in\mathcal X$, the process $Z^x$ is stochastically
continuous.

The corresponding transition semigroup is given by $P_t\psi(x)=
\mathbb{E}
\psi(Z^x(t))$, $\psi\in B_b(\mathcal X)$, and we assume
that it is Feller.
\begin{definition}
$\Phi\dvtx\mathcal X\to[0,+\infty)$ is called a \textit{Lyapunov
function} if it is measurable, bounded on bounded sets and
$
\lim_{\|x\|_{\mathcal X}\uparrow\infty}\Phi(x)=\infty.
$
\end{definition}

We shall assume that the deterministic equation
%
%
\begin{equation}\label{082102}
\frac{\dd Y(t)}{\dd t} = AY(t) + F(Y(t)),\qquad Y(0)=x,
\end{equation}
defines a continuous semi-dynamical system $(Y^x, x\in\mathcal X)$, that
is, for each \mbox{$x\in\mathcal X$}, there
exists a unique continuous solution to (\ref{082102}) that we denote
by $Y^x=(Y^x(t))_{ t\ge0}$ and for a given $t$, the mapping
$x\mapsto Y^x(t)$ is measurable. Furthermore, we have
$Y^{Y^x(t)}(s)=Y^x(t+s)$ for all $t,s\ge0$ and $x\in\mathcal X$.
\begin{definition}
\label{def4.2} A set ${\mathcal K}\subset\mathcal X$ is called \textit{a
global attractor} for the semi-dynamical system if:
\begin{enumerate}[(2)]
\item[(1)] it is invariant under the semi-dynamical system, that is,
$Y^x(t)\in\mathcal K$ for any $x\in\mathcal K$ and $t\ge0$;
\item[(2)] for any $\varepsilon,R>0$, there exists
$T$ such that $Y^{x}(t)\in\mathcal K+\varepsilon B(0,1)$ for $t\ge T$ and
$\|x\|_{\mathcal X}\le R$.
\end{enumerate}
\end{definition}
\begin{definition}
\label{def4.3} The family $(Z^x(t))_{t\ge0}, x\in\mathcal X$, is
\textit{stochastically stable} if
%
%
\begin{equation}\label{ET240}
\forall \varepsilon,  R,  t>0\qquad \inf_{x\in
B(0,R)}\mathbb{P} \bigl(\|Z^x(t)-Y^x (t)\|_{\mathcal X}<\varepsilon
\bigr)>0.
\end{equation}
\end{definition}

In Section \ref{SPDEs}, from Theorems \ref{Mtheorem} and
\ref{Theorem}, we derive the following result concerning ergodicity of $Z$.
\begin{theorem}\label{TGeneral}
Assume that:
\begin{longlist}
\item the semi-dynamical system $(Y^x, x\in\mathcal X)$ defined by
(\ref{082102})
has a compact, global attractor $\mathcal K$;
\item
$(Z^x(t))_{t\ge0}$ admits a Lyapunov function $\Phi$, that is,
\[
\forall x\in
\mathcal X\qquad \sup_{t\ge0}\mathbb E\Phi(Z^x (t)) <\infty;
\]
\item the family $(Z^x(t))_{t\ge0}, x\in\mathcal X$, is
stochastically stable and
%
%
\begin{equation}\label{condiii}
\bigcap_{x\in\mathcal K}\bigcup_{t\ge0} \Gamma^t(x)\neq\varnothing,
\end{equation}
where $\Gamma^t(x)=\operatorname{supp}P_t^*\delta_{x}$; 
%
\item its transition semigroup has the e-property.
\end{longlist}
$(Z^x(t))_{t\ge0}, x\in\mathcal X$, then admits a unique, invariant
measure $\mu_*$
and is weak-$^*$ mean ergodic. Moreover, for any bounded, Lipschitz
observable $\psi$, the weak law of large numbers holds:
\[
\mathop{{\mathbb{P}} \mbox{-}\mathrm{lim}}_{T\to+\infty}\frac{1}{T}\int
_0^T\psi(Z^x(t))\,\dd t=\int_{\mathcal X}\psi\,\dd\mu_*.
\]
\end{theorem}
\begin{remark}
Observe that condition (\ref{condiii}) in Theorem
\ref{TGeneral} is trivially satisfied if $\mathcal K$ is a singleton. Also,
this condition holds if the semi-dynamical system,
obtained after removing the noise, admits a global attractor that is
contained in the support of the
transition probability function of the solutions of (\ref{E1})
corresponding to the starting point at the attractor (this situation
occurs, e.g., if the noise is nondegenerate).

Another situation when (\ref{condiii}) can be guaranteed occurs if we
assume (\ref{ET240}) and uniqueness of an invariant probability measure
for $(Y^x, x\in\mathcal X)$. From stochastic stability condition
(\ref{ET240}), it is clear that the support of such a measure is
contained in any
$ \bigcup_{t\ge0}\Gamma_t(x)$ for $x\in\mathcal K$. We do not know,
however, whether there exists an example
of a semi-dynamical system corresponding to (\ref{082102}) with a
nonsingle point attractor and such that it admits a unique invariant measure.
\end{remark}
\begin{remark}
The e-property used in Theorem
\ref{TGeneral} can be understood as an intermediary between the
strong dissipativity property of \cite{DZ2} and asymptotic strong
Feller property
(see \cite{HM}). A trivial
example of a transition probability semigroup that is neither
dissipative (in the sense of \cite{DZ2}) nor asymptotic strong
Feller, but satisfies the e-property,
is furnished by the dynamical system on a unit
circle $\{z\in\mathbb C\dvtx|z|=1\}$ given by $\dot z=i\alpha z$,
where $\alpha/(2\pi)$ is an irrational real.
For more examples of Markov processes
that have the e-property, but are neither
dissipative nor have the
asymptotic strong Feller property, see \cite{LasotaSzarek}.
A careful analysis of the current proof shows that the e-property could
be viewed as a consequence of a certain version of the asymptotic
strong Feller property concerning time averages of the transition
operators. We shall investigate this point in more detail in a
forthcoming paper.
\end{remark}

Our last result follows from an application of the above theorem and
concerns the weak law of large numbers for the passive tracer in a
compressible random flow. The trajectory of a particle is then described
by the solution of an ordinary differential equation,
%
%
\begin{equation}\label{ET1}
\frac{\dd\mathbf x(t)}{\dd t}=V(t,\mathbf x(t)),\qquad \mathbf
x(0)=\mathbf x_0,
\end{equation}
where $V(t,\xi), (t,\xi)\in\mathbb{R}^{d+1}$, is a $d$-dimensional
random vector field. This is a simple model used in statistical
hydrodynamics that describes transport of matter in a turbulent
flow. We assume that $V(t,\xi)$ is mean zero, stationary,
spatially periodic, Gaussian and Markov in a time random field. Its
covariance matrix
\[
R_{i,j}(t-s,\xi-\eta) :=\mathbb{E} [V_i(t,\xi)V_j(s,\eta)]
\]
is given by its Fourier coefficients,
\begin{eqnarray*}
\widehat R_{i,j}(h,k) :\!&=&\frac{1}{(2\pi)^d}\int_{\mathbb{T}^d}
e^{-\ii
k\xi}
R_{i,j}(h,\xi) \,\dd\xi\\
&=& e^{-\gamma(k)|h|
}\mathcal E_{i,j}(k),\qquad i,j=1,\ldots,d,
\end{eqnarray*}
$h\in\mathbb{R}, k\in\mathbb{Z}^d$. Here, $\mathbb{T}^d:=[0,2\pi
)^d$, the \textit
{energy spectrum}
$\mathcal E:=[\mathcal E_{i,j}]$ maps $\mathbb{Z}^d$ into the space
$S_+(d)$ of all
nonnegative definite Hermitian matrices and the \textit{mixing
rates} $\gamma\dvtx\mathbb{Z}^d\to(0,+\infty)$. Denote by
$\operatorname{Tr} A$ the trace of a given $d\times d$ matrix $A$
and by
${{\mathbb{P}} \mbox{-}\mathrm{lim}}$ the limit in probability. In Section
\ref{sec5}, we show the following result.
\begin{theorem}\label{TTracer}
Assume that
%
%
\begin{eqnarray}\label{H1}
&\displaystyle\exists m>d/2+1, \alpha\in(0,1)\qquad \nnorm\mathcal E\nnorm^2:=
\sum_{k\in\mathbb{Z}^d} \gamma^{\alpha} (k)|k|^{2(m+1)}
\operatorname{Tr}
\mathcal E(k)
<\infty,\hspace*{-37pt}&
\\
%
%
\label{H2}
&\displaystyle\int_0^\infty\sup_{k\in\mathbb{Z}^d} e^{-\gamma(k)t}
|k| \,\dd t
<\infty.&
\end{eqnarray}
There then exists a constant vector $v_*$ such that
\[
\mathop{{\mathbb{P}} \mbox{-}\mathrm{lim}}_{t\uparrow
\infty}\frac{\mathbf x(t)}{t}= v_*.
\]
\end{theorem}
\begin{remark}
We will show that $v_*=\mathbb{E}_{\mu_*} V(0,0)$, where the
expectation $\mathbb{E}_{\mu_*}$ is calculated with respect to the path
measure that corresponds to the Markov process starting with the
initial distribution $\mu_*$, which is invariant under Lagrangian
observations of the velocity field, that is, the vector field-valued
process $V(t,\mathbf x(t)+\cdot)$, $t\ge0$. In the physics literature,
$v_*$ is referred to as the \textit{Stokes drift}. Since $V$ is
spatially stationary, the Stokes drift does not depend on the
initial value $\mathbf x_0$.
\end{remark}
\begin{remark}
Note that condition (\ref{H2}) holds if
\[
\exists \varepsilon,K_0>0 \mbox{ such that } \forall k\in
\mathbb{Z}
^d_*\qquad \gamma(k)\ge K_0|k|^{1+\varepsilon}.
\]
Indeed, it is clear that, under this assumption,
\[
\int_1^\infty\sup_{k\in\mathbb{Z}^d_*} e^{-\gamma(k)t} |k|
\,\dd
t <\infty.
\]
On the other hand, for $t\in(0,1]$, we obtain
\[
\sup_{k\in\mathbb{Z}^d_*} e^{-\gamma(k)t} |k|\le\sup_{k\in
\mathbb{Z}^d_*}
\exp\{-K_0|k|^{1+\varepsilon}t+\log|k| \}\le\frac
{C}{t^{1/(1+\varepsilon)}}
\]
for some constant $C>0$. This, of course, implies (\ref{H2}).
\end{remark}

\section[Proofs of Theorems 1 and 2]{Proofs of Theorems \protect\ref{Mtheorem} and \protect\ref
{Theorem}}

\subsection{Some auxiliary results}

For the proof of the following lemma the reader is referred to \cite
{LasotaSzarek};
see the argument given on pages
517 and 518.
\begin{lemma}
\label{szarek-lasota} Suppose that $(\nu_n)\subset{\mathcal M}_1$
is not tight. There then exist an $\varepsilon>0$, a sequence of compact
sets $(K_i)$ and an increasing sequence of positive integers
$(n_i)$ satisfying
%
%
\begin{equation}\label{A1g}
\nu_{n_i}(K_i)\ge\varepsilon \qquad\forall i,
\end{equation}
and
%
%
\begin{equation}\label{A2}
\min\{\rho(x, y)\dvtx x\in K_i, y\in K_j\}\ge\varepsilon
\qquad\forall i\neq j.
\end{equation}
\end{lemma}

Recall that $\mathcal{T}$ is defined by (\ref{Pr1}).
\begin{proposition}\label{Prop}
Suppose that $(P_t)_{t\ge0}$ has the e-property and admits an
invariant probability
measure $\mu_*$. Then $\operatorname{supp}\mu_*\subset\mathcal
T$.
\end{proposition}
\begin{pf}
Let $\mu_*$ be the invariant measure in question. Assume,
contrary to our claim, that $ (Q^T(x) )_{ T\ge0}$ is not
tight for some $x\in\operatorname{supp}\mu_*$. Then, according to Lemma
\ref{szarek-lasota}, there exist a strictly increasing sequence of
positive numbers $T_i\uparrow\infty$, a~positive number $\varepsilon
$ and a
sequence of compact sets $(K_i)$ such that
%
%
\begin{equation}\label{A1}
Q^{T_i}(x,K_i)\ge\varepsilon\qquad\forall i,
\end{equation}
and (\ref{A2}) holds.
We will derive the assertion from the claim that there exist
sequences $(\tilde{f}_n)\subset\operatorname{Lip}_b (\mathcal X)$,
$(\nu_n)\subset{\mathcal M}_1$ and an increasing sequence of
integers $(m_n)$ such that $\operatorname{supp}\nu_n\subset B(x,
1/n)$ for any
$n$, and
%
%
\begin{equation}\label{A3}
\mathbf{1}_{K_{m_n}}\le\tilde{f}_n\le
\mathbf{1}_{K_{m_n}^{\varepsilon/4}} \quad\mbox{and}\quad
\operatorname{Lip} (\tilde{f}_n)\le4/\varepsilon\qquad\forall n.
\end{equation}
Here, $A^\varepsilon:=\{x\in\mathcal X\dvtx \operatorname{dist}
(x,A)<\varepsilon\}$, with
$\varepsilon>0$, denotes the $\varepsilon$-neighborhood of $A\subset
\mathcal X$. Moreover,
%
%
\begin{equation}\label{A4}
{P}_t^*\nu_n \Biggl(\bigcup_{i=n}^{\infty}
K_{m_i}^{\varepsilon/4} \Biggr)\le\varepsilon/4 \qquad\forall
t\ge0,
\end{equation}
and
%
%
\begin{equation}\label{A5}
|{P}_t f_n (x)-{P}_t f_n (y)|<\varepsilon/4 \qquad\forall t\ge0,
\forall y\in\operatorname{supp}\nu_n,
\end{equation}
$f_1:=0$ and $f_n:={\sum_{i=1}^{n-1}} \tilde{f}_i$, $n\ge2$.
Temporarily admitting the above claim, we show how to complete the
proof of the
proposition. First, observe that (\ref{A2}) and condition (\ref{A3})
together imply that the series $f:=\sum_{i=1}^{\infty} \tilde{f}_i$
is uniformly convergent and $\|f\|_\infty=1$. Also, note that for
$x,y$ such that $\rho(x,y)<\varepsilon/8$, we have \mbox{$\tilde f_i(x)\not
=0$}, or
$\tilde f_i(y)\not=0$, for at most one $i$. Therefore, for such points,
$|f(x)-f(y)|<16\varepsilon^{-1}\rho(x,y)$. This, in particular,
implies that
$f\in\operatorname{Lip} (\mathcal X)$. From (\ref{A1}) and (\ref
{A3})--(\ref{A5}), it
follows that
%
%
\begin{eqnarray}\label{Add5}
&&\int_{\mathcal X} Q^{T_{m_n}}(x,\dd y) f(y) - \int_{\mathcal X}
Q^{T_{m_n}}\nu_{n}(\dd y) f(y)\nonumber\\
&&\qquad \ge Q^{T_{m_n}}(x,K_{m_n})
+\int_{\mathcal X} Q^{T_{m_n}}(x,\dd y) f_{n}(y)\\
&&\qquad\quad{} -\int_\mathcal X
Q^{T_{m_n}}\nu_n(\dd y) f_{n}(y)-Q^{T_{m_n}}\nu_n
\Biggl(\bigcup_{i=n}^{\infty} K_{m_i}^{\varepsilon/4} \Biggr).\nonumber
\end{eqnarray}
By virtue of (\ref{A1}), the first term on the right-hand side of
(\ref{Add5}) is greater than or equal to $\varepsilon$. Combining the second
and the third terms, we obtain that their absolute value equals
\[
\biggl|\frac{1}{T_{m_n}}\int_0^{T_{m_n}}\int_\mathcal X [P^sf_{n}(x)-
P^sf_{n}(y) ]\nu_n(\dd y)\,\dd
s \biggr|\stackrel{(\fontsize{8.36}{10}\selectfont{\mbox{\ref{A5}}})}{\le}\frac{\varepsilon}{4}.
\]
The fourth term is less than or equal to $\varepsilon/4$, by virtue of
(\ref{A4}).
Summarizing, we have shown that
\begin{eqnarray*}
&&
\int_\mathcal XQ^{T_{m_n}}(x,\dd y) f(y) - \int_\mathcal
XQ^{T_{m_n}}\nu_{n}(\dd
y) f(y)
\\
&&\qquad= \frac{1}{T_{m_n}}\int_0^{T_{m_n}}\dd s\int_\mathcal X [P^sf(x)-
P^sf(y) ]\nu_n(\dd y) >\frac{\varepsilon}{2}
\end{eqnarray*}
for every positive integer $n$. Hence, there must be a sequence
$(t_n,y_n)$ such that $t_n\in[0,T_{m_n}]$, $y_n\in\operatorname{supp}
\nu_n\subset B(x,1/n)$, for which $
{P}_{t_n}f(x)-{P}_{t_n}f(y_n)>\varepsilon/2, $ $n\ge1$. This clearly
contradicts equicontinuity of $(P_t f)_{t\ge0}$ at $x$.
\begin{pf*}{Proof of the claim}
We accomplish this by induction on $n$. Let
$n=1$. Since $x\in\operatorname{supp}\mu_*$, we have $\mu_* (B(x,
\delta))>0$ for
all $\delta>0$. Define the probability measure $\nu_1$ by the formula
\[
\nu_1 (B)=\mu_*(B | B(x, 1)) := \frac{\mu_*(B\cap
B(x,1))}{\mu_*(B(x,1))},\qquad B\in\mathcal B (\mathcal X).
\]
Since $\nu_1\le\mu_*^{-1}( B(x, 1))\mu_*$, from the fact that
$\mu_*$ is invariant, it follows that the family $(P_t^*\nu_1)_{t\ge
0}$ is tight.
Thus, there exists a compact set $K$ such that
%
%
\begin{equation}\label{082101}
{P}_t^*\nu_1 (K^c)\le\varepsilon/4\qquad \forall t\ge0.
\end{equation}
Note, however, that
$K\cap K_i^{\varepsilon/4}\not=\varnothing$ for only finitely many $i$'s.
Otherwise, in light of~(\ref{A2}), one
could construct in $K$ an infinite set of points separated from each other
at a distance of at least $\varepsilon/2$, which contradicts its compactness.
As a result, there exists
an integer $m_1$ such that
\[
{P}_t^*\nu_1 \Biggl(\bigcup_{i=m_1}^{\infty}
K_i^{\varepsilon/4} \Biggr)\le\varepsilon/4 \qquad\forall t\ge0.
\]
Let\vspace*{-2pt} $\tilde{f}_1$ be an arbitrary Lipschitz function satisfying $
\mathbf{1}_{K_{m_1}}\le\tilde{f}_1\le\mathbf{1}_{K_{m_1}^{\varepsilon/4}}
$ and $\operatorname{Lip} (\tilde{f}_1)\le4/\varepsilon$.

Assume, now, that for a given $n\ge1$, we have already constructed
$\tilde{f}_1,\ldots, \tilde{f}_{n}$, $\nu_1,\ldots, \nu_{n}$,
$m_1,\ldots, m_{n}$ satisfying (\ref{A3})--(\ref{A5}). Since $(P_t
f_{n+1})_{t\ge0}$ is equicontinuous, we can choose $\delta<1/(n+1)$
such that $|{P}_t f_{n+1} (x)-{P}_t f_{n+1} (y)|<\varepsilon/4$ for
all $t\ge0$ and $y\in B(x,\delta)$. Suppose, further, that
$\nu_{n+1}:=\mu_*(\cdot| B(x, \delta))$. Since the measure is
supported in $B(x, \delta)$, condition (\ref{A5}) holds for $f_{n+1}$.
Tightness of $(P_t^*\nu_{n+1})_{t\ge0}$ can be argued in the same
way as in the case $n=1$. As a consequence, one can find $m_{n+1}>m_n$
such that
\[
{P}_t^*\nu_{n+1} \Biggl(\bigcup_{i=m_{n+1}}^{\infty}K_i^{\varepsilon
/4} \Biggr)
\le\varepsilon/4 \qquad\forall t\ge0.
\]
Finally, we let $\tilde{f}_{n+1}$ be an arbitrary continuous
function satisfying (\ref{A3}).
\end{pf*}
\noqed\end{pf}

For given an integer $k\ge1$, times $t_1,\ldots, t_k\ge0$ and
a measure $\mu\in\mathcal M_1$, we let $ Q^{t_k,\ldots, t_1}\mu
:=Q^{t_k}\cdots Q^{t_1}\mu$. The following simple lemma will be
useful in the sequel. In what follows, $\|\cdot\|_{\mathrm{TV}}$
denotes the total variation norm.
\begin{lemma}\label{Lemma}
For all $k\ge1$ and $t_1,\ldots, t_k>0$,
%
%
\begin{equation}\label{Le1}
{\limsup_{T\to+\infty}\sup_{\mu\in\mathcal M_1}}
\|Q^{T,t_k,\ldots, t_1}\mu-Q^T\mu\|_{\mathrm{TV}}=0.
\end{equation}
\end{lemma}
\begin{pf} To simplify the notation, we assume that $k=1$.
The general case can be argued by the induction on the length of
the sequence $t_1, \ldots, t_k$ and is left to the reader. For any
$T>0$, we have
\begin{eqnarray*}
Q^{T,t_1}\mu-Q^T\mu
&=&
{(T t_1)}^{-1}\int_0^{t_1}\dd r\,
\biggl[\int_{0}^{T} {P}_{s+r}^*\mu\,\dd s -\int_{0}^{T} {P}_{s}^*\mu\,\dd
s \biggr] \\
&=&
{(Tt_1)}^{-1}\int_0^{t_1}\dd r\int_{0}^{r} ({P}_{s+T}^*\mu-P_s^*\mu
)\,\dd s .
\end{eqnarray*}
%
The total variation norm of $Q^{T,t_1}\mu-Q^T\mu$ can therefore be
estimated by
$t_1/T$ and (\ref{Le1}) follows.
\end{pf}
%


\subsection[Proof of Theorem 1]{Proof of Theorem \protect\ref{Mtheorem}}
\label{sec2.3} The existence of an invariant measure follows from
Theorem 3.1 of \cite{LasotaSzarek}. We will show that for arbitrary
$x_1, x_2\in\mathcal T$ and $\psi\in\operatorname{Lip}_b (X)$,
%
%
\begin{equation}\label{U1}
\lim_{T\uparrow\infty} \biggl|\int_\mathcal X\psi
(y)Q^T(x_1,\dd
y)-\int_\mathcal X\psi(y)Q^T(x_2,\dd y) \biggr|=0 .
\end{equation}
From this, we can easily deduce (\ref{Mt2}) using, for instance,
Example 22,
page 74 of~\cite{pollard}. Indeed, for any $\nu$, as in the statement
of the theorem,
\begin{eqnarray*}
&&\int_\mathcal X\psi(y)Q^T\nu(\dd y)-\int_\mathcal X\psi\,\dd\mu_*
\\
&&\qquad= \int_\mathcal X\int_\mathcal X\nu(\dd x) \mu_*(\dd x') \biggl(\int
_\mathcal X
\psi(y)Q^T(x,\dd y)-\int_\mathcal X\psi(y)Q^T(x',\dd y) \biggr)
\end{eqnarray*}
and (\ref{Mt2}) follows directly from (\ref{U1}) and Proposition \ref{Prop}.
The rest of the argument will be
devoted to the proof of (\ref{U1}).

Fix a sequence $(\eta_n)$ of positive numbers monotonically
decreasing to $0$. Also, fix arbitrary $\varepsilon>0$, $\psi\in
\operatorname{Lip}_b (\mathcal X)$, $x_1,x_2\in\mathcal T$. For
these parameters,
we define $\Delta\subset\mathbb R$ in the following way:
$\alpha\in\Delta$ if and only if $\alpha>0$ and there exist a
positive integer~$N$, a sequence of times $(T_{\alpha, n})$
and sequences of measures $(\mu_{\alpha, i}^n), (\nu_{\alpha,
i}^n)\subset\mathcal M_1$, $i=1, 2$, such that for $n\ge N$,
%
%
%
%
\begin{eqnarray}\label{U1b}
T_{\alpha, n}&\ge& n,
\\
\label{U2}
\|Q^{T_{\alpha, n}}(x_i)-\mu_{\alpha, i}^n\|_{\mathrm{TV}}&<&\eta_n,
\\
%
%
\label{U3}
\mu_{\alpha, i}^n&\ge&\alpha\nu_{\alpha, i}^n \qquad\mbox
{for $i=1, 2$},
\end{eqnarray}
and
%
%
\begin{equation}\label{U4}
\limsup_{T\uparrow\infty} \biggl| \int_\mathcal X
\psi(x)Q^{T}\nu_{\alpha, 1}^n(\dd x) - \int_\mathcal X
\psi(x)Q^{T}\nu_{\alpha, 2}^n(\dd x) \biggr|<\varepsilon.
\end{equation}
Our main tool is contained in the following lemma.
\begin{lemma}
\label{delta} For given $\varepsilon>0$, $(\eta_n)$, $x_1,x_2\in
\mathcal T$ and $\psi\in\operatorname{Lip}_b (\mathcal X)$, the set
\mbox{$\Delta\not=\varnothing$}. Moreover, we have
$\sup\Delta=1$.
\end{lemma}

Accepting the truth of this lemma, we show how to complete the proof of
(\ref{U1}). To that end, let us choose an arbitrary $\varepsilon>0$.
Then there exists an $\alpha>1-\varepsilon$ that belongs to $\Delta
$. By
virtue of (\ref{U2}), we can replace the $Q^T(x_i,\cdot)$ appearing
in (\ref{U1})
by $\mu_{\alpha, i}^n$ and the resulting error can be estimated
for
$T\ge T_{\alpha,n}$ as follows:
%
%
\begin{eqnarray}\label{073105}\qquad
&& \biggl| \int_\mathcal X\psi(y)Q^{T}(x_1,\dd y) -
\int_\mathcal X\psi(y)Q^{T}(x_2,\dd y) \biggr|
\nonumber\\
&&\qquad \le\sum_{i=1}^2 \biggl| \int_\mathcal X\psi
(y)Q^{T}(x_i,\dd
y) -
\int_\mathcal X\psi(y)Q^{T,T_{\alpha,n}}(x_i,\dd y) \biggr|
\nonumber\\
&&\qquad\quad{} + \biggl| \int_\mathcal X\psi(y)Q^{T}\mu_{\alpha
,1}^n(\dd y) -
\int_\mathcal X\psi(y)Q^{T}\mu_{\alpha,2}^n(\dd y) \biggr|
\nonumber\\
&&\qquad\quad{} +\sum_{i=1}^2 \biggl| \int_\mathcal X
\psi(y)Q^{T,T_{\alpha,n}}(x_i,\dd y) -
\int_\mathcal X\psi(y)Q^{T}\mu_{\alpha,i}^n(\dd y) \biggr|
\\
&&\qquad \le\sum_{i=1}^2 \biggl| \int_\mathcal X\psi
(y)Q^{T}(x_i,\dd
y) -
\int_\mathcal X\psi(y)Q^{T,T_{\alpha,n}}(x_i,\dd y) \biggr|
\nonumber\\
&&\qquad\quad{}
+ \biggl| \int_\mathcal X\psi(y)Q^{T}\mu_{\alpha,1}^n(\dd y) -
\int_\mathcal X\psi(y)Q^{T}\mu_{\alpha,2}^n(\dd y) \biggr|\nonumber\\
&&\qquad\quad{}+2\eta_n\|
\psi
\|_\infty.\nonumber
\end{eqnarray}
To deal with the second term on the last right-hand side of
(\ref{073105}), we use condition (\ref{U3}). We can then replace
$\mu_{\alpha,i}^n$ by $\nu_{\alpha,i}^n$ and obtain
%
%
\begin{eqnarray}\label{073106}\quad
&& \biggl| \int_\mathcal X\psi(y)Q^{T}\mu_{\alpha,1}^n(\dd y) -
\int_\mathcal X\psi(y)Q^{T}\mu_{\alpha,2}^n(\dd y) \biggr|
\nonumber\\
&&\qquad\stackrel{(\fontsize{8.36}{10}\selectfont{\mbox{\ref{U3}}})}{\le} \alpha\biggl| \int_\mathcal X
\psi(y)Q^{T}\nu_{\alpha,1}^n(\dd y) - \int_\mathcal X
\psi(y)Q^{T}\nu_{\alpha,2}^n(\dd
y) \biggr|\nonumber\\[-8pt]\\[-8pt]
&&\qquad\quad\hspace*{6.5pt}{} + \sum_{i=1}^2\|\psi\|_{\infty}
(\mu_{\alpha,i}^n-\alpha\nu_{\alpha,i}^n)(\mathcal X)
\nonumber\\
&&\hspace*{6.5pt}\qquad\le\biggl| \int_\mathcal X\psi(y)Q^{T}\nu_{\alpha,1}^n(\dd y) -
\int_\mathcal X\psi(y)Q^{T}\nu_{\alpha,2}^n(\dd y)
\biggr|+2\varepsilon\|\psi\|
_{\infty}.\nonumber
\end{eqnarray}
In the last inequality, we have used the fact that $1-\alpha
<\varepsilon$.
Summarizing, from Lemma \ref{Lemma}, (\ref{073105}), (\ref{073106})
and (\ref{U4}), we obtain that
\[
\limsup_{T\uparrow\infty} \biggl| \int_\mathcal X
\psi(y)Q^{T}(x_1,\dd y) -
\int_\mathcal X\psi(y)Q^{T}(x_2,\dd y) \biggr|\le
2\eta_n\|\psi\|_{\infty}+2\varepsilon\|\psi\|_{\infty
}+\varepsilon.
\]
Since $\varepsilon>0$ and $n$ were arbitrarily chosen, we conclude
that (\ref{U1}) follows.
\begin{pf*}{Proof of Lemma \protect\ref{delta}}
First, we show that $\Delta\ne\varnothing$. Let $z\in\mathcal X$ be
such that
for every $\delta>0$ and $x\in\mathcal X$, condition (\ref{Mt1}) is
satisfied. Equicontinuity of $(P_t\psi)_{t\ge0}$ at $z\in\mathcal X$ implies
the existence of $\sigma>0$ such that
%
%
\begin{equation}
\label{U4b}
|{P}_{t} \psi(z)-{P}_t \psi(y)|<\varepsilon/2 \qquad\mbox{for $y\in
B(z,\sigma)$ and $t\ge0$.}
\end{equation}
By (\ref{Mt1}), there exist $\beta>0$ and $T_0>0$ such that
%
%
\begin{equation}
\label{072701}
Q^T(x_i, B(z,\sigma))\ge\beta\qquad\forall T\ge T_0, i=1, 2.
\end{equation}
Set $\alpha:=\beta$ and $T_{\alpha, n}=n+T_0$ for $n\in\mathbb
{N}$, $
\mu_{\alpha, i}^n:=Q^{T_{\alpha,n}}(x_i) $ and $ \nu_{\alpha, i}^n
(\cdot):=\mu_{\alpha, i}^n (\cdot| B(z,\sigma)) $ for $i=1, 2$ and
$n\ge1$. Note that $\mu_{\alpha, i}^n ( B(z,\sigma))>0$, thanks
to~(\ref{072701}). The measures $\nu_{\alpha, i}^n$, $i=1,2$, are supported
in $ B(z,\sigma)$ and, therefore, for all $t\ge0$,
we have
\begin{eqnarray*}
&& \biggl|\int_\mathcal X\psi(x) P_t^*\nu_{\alpha, 1}^n(\dd x)-
\int_\mathcal X\psi(x) P_t^*\nu_{\alpha, 2}^n(\dd x) \biggr|\\
&&\qquad =
\biggl|\int_\mathcal XP_t \psi(x) \nu_{\alpha, 1}^n(\dd x)-
\int_\mathcal XP_t\psi(x) \nu_{\alpha, 2}^n(\dd x) \biggr|\\
&&\qquad \le\biggl|\int_\mathcal X [P_t\psi(x)-P_t\psi(z) ]
\nu_{\alpha, 1}^n(\dd x) \biggr| \\
&&\qquad\quad{} + \biggl|\int_\mathcal X [P_t
\psi(x)-P_t\psi(z) ] \nu_{\alpha, 2}^n(\dd
x) \biggr|\stackrel{(\fontsize{8.36}{10}\selectfont{\mbox{\ref{U4b}}})}{<}\varepsilon.
\end{eqnarray*}
Hence, (\ref{U4}) follows. Clearly, conditions (\ref{U1})--(\ref
{U3}) are
also satisfied. Thus, $\Delta\ne\varnothing$.

Next, we show that $ \sup\Delta=1. $ Suppose, contrary to
our claim, that $\alpha_0:=\sup\Delta<1$. Thanks to the previous
step, we have $\alpha_0>0$. Let $(\alpha_n)\subset\Delta$ be such that
$\lim_{n\to\infty}\alpha_n=\alpha_0$. Set $ T_n:=T_{\alpha_n, n}$,
$\mu_{n, i}:=\mu_{\alpha_n, i}^n$ and $\nu_{n, i}:=\nu_{\alpha_n, i}^n$
for $n\ge1$ and $i=1, 2$. From conditions (\ref{U2}), (\ref{U3}) and the
fact that the family $(Q^T(x_i))$ is tight for $i=1, 2$, it
follows that the sequences $(\mu_{n, i})$, $(\nu_{n, i})$, $i=1, 2$,
are also tight. Indeed, (\ref{U2}) clearly implies tightness of
$(\mu_{n, i})$, $i=1, 2$. As a consequence, for any $\varrho>0$, there
exists a compact set $K\subset\mathcal X$ such that $\mu_{n,
i}(\mathcal X\setminus K)<\varrho$ for all $n\ge1$, $i=1,2$. In turn,
condition
(\ref{U3}) implies that for sufficiently large $n$, we have
\[
\nu_{n, i}(\mathcal X\setminus K)< \frac{2\mu_{n, i}(\mathcal
X\setminus
K)}{\alpha_0}<\frac{2\varrho}{\alpha_0}
\]
and tightness of $(\nu_{n, i})$, $i=1, 2$, follows. Therefore,
without loss of generality, we may assume that the sequences
$(\mu_{n, i})$, $(\nu_{n, i})$, $i=1, 2$, are weakly convergent.
The sequences
%
%
\begin{equation}\label{072803}
\bar\mu_{n,i}:=\mu_{n, i}-\alpha_n\nu_{n, i},\qquad
n\ge1,
\end{equation}
are therefore also weakly convergent for $i=1,2$. The assumption that
$\alpha_0<1$ implies that the respective limits are nonzero
measures; we denote them by $\bar\mu_i$, $i=1,2$, correspondingly.
Let $y_i\in\operatorname{supp}\bar\mu_i$, $i=1,2$. Analogously to
the previous
step, we may choose $\sigma>0$ such that (\ref{U4b}) is satisfied.
By (\ref{Mt1}), we choose $T>0$ and $\gamma>0$ for which
%
%
\begin{equation}\label{072704}
Q^T\bigl(y_i,B(z, \sigma/2)\bigr)\ge\gamma\qquad\mbox{for $i=1, 2$}.
\end{equation}
Since the semigroup $(P_t)_{t\ge0}$ is Feller,
we may find $r>0$ such that
%
%
\begin{equation}\label{072705}
Q^T(y,B(z, \sigma))\ge\gamma/2 \qquad\mbox{for $y\in B(y_i, r)$ and
$i=1, 2$}.
\end{equation}
Indeed, it suffices to choose $\phi\in\operatorname{Lip}_b
(\mathcal X)$ such
that $\mathbf{1}_{B(z, \sigma/2)}\le\phi\le\mathbf{1}_{B(z,
\sigma)}$. From
(\ref{072704}), we have $\int_\mathcal X\phi(x) Q^T(y_i,\dd x)\ge
\gamma$.
The Feller property implies that there exists $r>0$ such that, for
$y\in
B(y_i, r)$ and $i=1, 2$, we have
\[
Q^T(y,B(z, \sigma))\ge\int_{\mathcal X} \phi( x) Q^T(y,\dd x)\ge
\frac
\gamma2.
\]
Set
$s_0=\min\{\bar\mu_1 (B(y_1, r)), \bar\mu_2 (B(y_2, r))\}>0$.
Using part (iv) of Theorem 2.1, page 16 of \cite{billingsley}, we may
find $N\ge1$ such that
%
%
\begin{equation}
\label{080101} \bar\mu_{n, i} (B(y_i, r))>\frac{s_0}2 \quad\mbox{and}\quad
\alpha_n+s_0\frac\gamma4>\alpha_0
\end{equation}
for $n\ge N$. We prove that $\alpha_0':=\alpha_0+s_0\gamma/8$ also belongs
to $\Delta$, which obviously leads to a contradiction with the
hypothesis that $\alpha_0=\sup\Delta$. We construct sequences
$(T_{\alpha_0', n})$, $( \mu_{\alpha_0', i}^n)$ and
$(\nu_{\alpha_0', i}^n)$, $i=1,2$, that satisfy conditions
(\ref{U1b})--(\ref{U4}) with $\alpha$ replaced by $\alpha_0'$. Let
$\widehat\mu_n^i(\cdot):=\bar\mu_{n, i}(\cdot|B(y_i, r))$, $i=1,2$, be
the measure $\bar\mu_{n, i}$ conditioned on the respective balls
$B(y_i, r)$, $i=1,2$. That is, if $\bar\mu_{n, i}( B(y_i, r))\not
=0$, then we
let
%
%
\begin{equation}\label{U5}
\widehat\mu_n^i(\cdot) :=\frac{\bar\mu_{n, i}(\cdot\cap B(y_i, r))}
{\bar\mu_{n, i}( B(y_i, r))},
\end{equation}
while if $\bar\mu_{n, i}( B(y_i, r))=0$, we just let $\widehat
\mu_n^i(\cdot):=\delta_{y_i}$. Also, let $ \tilde\mu_n^i(\cdot):=(Q^T
\bar\mu_{n,i})\times\break (\cdot| B(z, \sigma))$.
%
From the above definition, it follows that
%
%
\begin{equation}\label{U5b}
Q^T\mu_{n,i}\ge\frac{s_0\gamma}{4}\tilde\mu_n^i+\alpha_n Q^T \nu
_{n, i}
\end{equation}
for $n\ge N$ and $i=1, 2$.
Indeed, note that from (\ref{080101}) and (\ref{U5}), we have
%
%
\begin{equation}
\label{072801}
\bar\mu_{n,i}(B)\ge
\frac{s_0}{2} \widehat\mu_n^i(B)\qquad \forall B\in{\mathcal
B}(\mathcal X),
\end{equation}
hence, also,
%
%
\begin{equation}\label{072802}
Q^T \bar\mu_{n,i}(B)\ge
\frac{s_0}{2} Q^T\widehat\mu_n^i(B) \qquad\forall B\in{\mathcal
B}(\mathcal X).
\end{equation}
On the other hand, by Fubini's theorem, we obtain
\begin{eqnarray*}
Q^T \widehat\mu_n^i (B(z, \sigma))
&=& T^{-1}\int_0^T\int_\mathcal X\mathbf{1}_{B(z, \sigma)}(x)
P_s^*\widehat\mu
_n^i(\dd x)\,\dd s\\
&=& T^{-1}\int_0^T\int_\mathcal XP_s\mathbf{1}_{B(z, \sigma)}(x)
\widehat\mu
_n^i(\dd
x)\,\dd s\\
&=&\int_\mathcal XQ^T(x,B(z, \sigma))\widehat\mu_n^i(\dd x)\\
&\stackrel{(\fontsize{8.36}{10}\selectfont{\mbox{\ref{U5}}})}{=}&\int_{B(y_i, r)}Q^T (x,B(z, \sigma))
\widehat\mu_n^i(\dd x) \stackrel{(\fontsize{8.36}{10}\selectfont{\mbox{\ref{072705}}})}{\ge} \frac{\gamma}{2}
\end{eqnarray*}
and, consequently, (\ref{072802}) implies that
%
%
\begin{equation}\label{072804}
Q^T \bar\mu_{n,i} (B(z, \sigma))\ge\frac{s_0\gamma}{4}.
\end{equation}
Hence, for any $B\in{\mathcal B}(\mathcal X)$,
\begin{eqnarray*}
Q^T\mu_{n,i}(B)
&\stackrel{(\fontsize{8.36}{10}\selectfont{\mbox{\ref{072803}}})}{=}&
Q^T\bar\mu_{n,i}(B)+\alpha_nQ^T\nu_{n,i}(B)\\
&\ge&
Q^T\bar\mu_{n,i}\bigl(B\cap B(z, \sigma)\bigr)+\alpha_nQ^T\nu_{n,i}(B)
\\
&\stackrel{(\fontsize{8.36}{10}\selectfont{\mbox{\ref{072804}}})}{\ge}&
\frac{s_0\gamma}{4} \tilde\mu_{n,i}(B)+\alpha_n Q^T\nu_{n,i}(B)
\end{eqnarray*}
and (\ref{U5b}) follows. At this point, observe that, by virtue of
(\ref{U5b}), measures $Q^T\mu_{n,i}$ and
$(s_0\gamma/4+\alpha_n)^{-1}[(s_0\gamma/4)\tilde\mu_{n,i}+\alpha_n
Q^T\nu_{n,i}]$ would satisfy (\ref{U3}), with $\alpha_0'$ in place of
$\alpha$, admitted them instead of $\mu_{\alpha_0',i}^n$ and
$\nu_{\alpha_0',i}^n$, respectively. Condition (\ref{U2}) need not,
however, hold in such case. To remedy this, we average
$Q^T\mu_{n,i}$ over a long time, using the operator $Q^R$ corresponding
to a sufficiently large $R>0$, and use Lemma \ref{Lemma}. More
precisely, since $\eta_n> \|Q^{T_n}(x_i)-\mu_{n,i}\|_{\mathrm{TV}}$ [thus,
also, $\eta_n> \|Q^{R,T,T_n}(x_i)-Q^{R,T}\mu_{n,i}\|_{\mathrm{TV}}$ for any
$R>0$], by Lemma \ref{Lemma}, we can choose $R_n>T_n$ such that
%
%
\begin{equation}\label{073101}\hspace*{28pt}
\|Q^{R_n,T,T_n}(x_i)-Q^{R_n}(x_i)\|_{\mathrm{TV}}
<\eta_n-
\|Q^{R_n,T,T_n}(x_i)-Q^{R_n,T}\mu_{n,i}\|_{\mathrm{TV}}.
\end{equation}
Let
%
%
\begin{equation} \label{U6}
T_{\alpha_0', n}:=R_n, \qquad \mu_{\alpha_0',
i}^n:=Q^{R_n}Q^T\mu_{n,i}
\end{equation}
and
%
%
\begin{equation}\label{U6b}
\nu_{\alpha_0', i}^n:=
\biggl(\alpha_n+\frac{s_0\gamma}{4} \biggr)^{-1}
Q^{R_n} \biggl(\alpha_nQ^T\nu_{n,i}+\frac{s_0\gamma}{4}\tilde\mu
_{n,i} \biggr)
\end{equation}
for $i=1, 2$, $n\ge1$.
By virtue of (\ref{073101}),
we immediately see that
\[
\|Q^{T_{\alpha_0', n}}(x_i)-\mu_{ \alpha_0', i}^n\|_{\mathrm{TV}}
<\eta_n \qquad\forall n\ge1.
\]
%
Furthermore, from (\ref{U5b}), positivity of $Q^{R_n}$ and the
definitions of $\alpha_0'$ and measures $\mu_{\alpha_0', i}^n$
$\nu_{\alpha_0', i}^n$, we obtain that
\[
\mu_{\alpha_0', i}^n\ge\alpha_0'\nu_{\alpha_0', i}^n
\qquad\forall
n\ge N, i=1, 2,
\]
when $N$ is chosen sufficiently large.
To verify (\ref{U4}), note that from (\ref{U6b}), it follows
that
%
%
\begin{eqnarray}\label{073102}
&&\biggl|\int_\mathcal X\psi(x)Q^S\nu_{\alpha_0', 1}^n (\dd x)-
\int_\mathcal X\psi(x)Q^S\nu_{\alpha_0', 2}^n (\dd x) \biggr|\nonumber\\
&&\qquad\le\alpha_n \biggl(\alpha_n+\frac{s_0\gamma}{4} \biggr)^{-1}\nonumber\\
&&\qquad\quad{}\times\biggl|\int_\mathcal X\psi(x) Q^{S,R_n,T}\nu_{n, 1}(\dd x)-
\int_\mathcal X\psi(x) Q^{S,R_n,T}\nu_{n, 2}(\dd x) \biggr| \\
&&\qquad\quad{}+\frac{s_0\gamma}{4} \biggl(\alpha_n+\frac{s_0\gamma}{4} \biggr)^{-1}
\biggl|\int_\mathcal X\psi(x) Q^{S,R_n}\tilde\mu_{n, 1}(\dd x)\,\dd s\nonumber\\
&&\qquad\quad\hspace*{101.2pt}{} -
\int_\mathcal X\psi(x)Q^{S,R_n}\tilde\mu_{n, 2}(\dd x) \biggr|\nonumber
\end{eqnarray}
for all $S\ge0$. Denote the integrals appearing in the first and
the second terms on the right-hand side of (\ref{073102}) by $I(S)$
and $\mathit{I I}(S)$, respectively. Condition (\ref{U4}) will follow if we
could demonstrate that the upper limits, as $S\uparrow\infty$, of
both of these terms are smaller than $\varepsilon$. To estimate
$I(S)$, we
use Lemma \ref{Lemma} and condition (\ref{U4}), which holds for
$\nu_{n, i}$, $i=1,2$. We then
obtain
\begin{eqnarray*}
\limsup_{S\uparrow\infty}I(S)&\le&\limsup
_{S\uparrow
\infty} \biggl|\int_\mathcal X\psi(x) Q^{S,R_n,T}\nu_{n, 1}(\dd x)-
\int_\mathcal X\psi(x) Q^{S}\nu_{n, 1}(\dd x) \biggr|
\\
&&{} +\limsup_{S\uparrow\infty} \biggl|\int_\mathcal X\psi(x)
Q^{S}\nu_{n, 1}(\dd x)- \int_\mathcal X\psi(x) Q^{S}\nu_{n, 2}(\dd
x) \biggr|
\\
&&{} +\limsup_{S\uparrow\infty} \biggl|\int_\mathcal X\psi(x)
Q^{S,R_n,T}\nu_{n, 2}(\dd x)- \int_\mathcal X\psi(x) Q^{S}\nu_{n,
2}(\dd
x) \biggr|<\varepsilon.
\end{eqnarray*}
%
On the other hand, since $\operatorname{supp}\tilde\mu_n^i\subset
B(z, \sigma)$,
$i=1,2$, we obtain, from equicontinuity condition (\ref{U4b}),
\begin{eqnarray*}
\mathit{I I}(S)
&=&\frac{1}{SR_n} \biggl|\int_0^S\int_0^{R_n}\int_\mathcal X\int
_\mathcal X
\bigl(P_{s_1+s_2}\psi(x)-P_{s_1+s_2}\psi(x')\bigr) \,ds_1\,ds_2\\
&&\hspace*{154.3pt}{}\times\tilde\mu_{n,
1}(\dd
x)\tilde\mu_{n, 2}(\dd x') \biggr|\le\frac{\varepsilon}{2}.
\end{eqnarray*}
Hence, (\ref{U4}) holds for $\nu_{\alpha'_0, i}^n$, $i=1,2$, and
function $\psi$. Summarizing, we have shown that
$\alpha_0'\in\Delta$. However, we also have
$\alpha_0'>\alpha_0=\sup\Delta$, which is clearly impossible.
Therefore, we conclude that $\sup\Delta=1$.
\end{pf*}


\subsection[Proof of Theorem 2]{Proof of Theorem \protect\ref{Theorem}}\label{sec2.4}

Taking Theorem \ref{Mtheorem} into account, the proof of the first
part of the
theorem will be completed as soon as we can show that $\mathcal
T={\mathcal X}$. Note that condition (\ref{Th1}) implies that $z\in$ supp
$\mu_*$. Indeed, let $B$ be a bounded set such that $\mu_*(B)>0$. We
can then write, for any $\delta>0$ and $T>0$,
\begin{eqnarray*}
\mu_*(B(z,\delta))&=&\int_\mathcal XQ^T(y,B(z,\delta))\mu
_*(\dd y)\\
&=&
\liminf_{T\uparrow\infty} \int_\mathcal X
Q^T(y,B(z,\delta))\mu_*(\dd y)
\\
&\stackrel{\mathrm{Fatou}\ \mathrm{lem.}}{\ge}& \int_\mathcal X
\liminf_{T\uparrow\infty}
Q^T(y,B(z,\delta))\mu_*(\dd y)\\
&\stackrel{(\fontsize{8.36}{10}\selectfont{\mbox{\ref{Th1}}})}{\ge}& \inf_{y\in B}
\liminf_{T\uparrow\infty} Q^T(y,B(z,\delta))\mu_*(B)>0.
\end{eqnarray*}
According to Proposition \ref{Prop}, the above implies that $z\in
{\mathcal T}$. Now, fix an arbitrary $x\in\mathcal X$.
Let $\mathcal C_{\varepsilon}$ be the family of all closed sets
$C\subset\mathcal X$ which possess a finite $\varepsilon$-net, that
is, there
exists a
finite set, say $\{x_1,\ldots, x_n\}$, for which
$C\subset\bigcup_{i=1}^n B(x_i, \varepsilon)$. To prove that the
family $(Q^T(x))$ is tight, it suffices to show that
for every $\varepsilon>0$, there exists $C_{\varepsilon}\in\mathcal
C_{\varepsilon}$ such that
%
%
\begin{equation}\label{Th3}
\liminf_{T\uparrow\infty}
Q^T(x,C_{\varepsilon})>1-\varepsilon;
\end{equation}
for more details, see, for example, pages 517 and 518 of \cite
{LasotaSzarek}. In
light of Lemma~\ref{Lemma}, this condition would follow if we could prove
that for given $\varepsilon>0$, $k\ge1$ and $t_1,\ldots,t_k\ge0$,
one can
find $T_\varepsilon>0$ and $C_{\varepsilon}\in\mathcal
C_{\varepsilon}$
such that
%
%
\begin{equation}\label{Th3b}
Q^{T,t_1,\ldots,t_k}(x,C_{\varepsilon})>1-\varepsilon\qquad\forall
T\ge T_\varepsilon.
\end{equation}

Fix an $\varepsilon>0$. Since $z\in{\mathcal T}$, we can find
$C_{\varepsilon/2}\in\mathcal C_{\varepsilon/2}$ such that (\ref{Th3})
holds with $\varepsilon/2$ in place of $\varepsilon$ and $x=z$. Let
$\tilde C:=C_{\varepsilon/2}^{\varepsilon/2}$ be the $\varepsilon
/2$-neighborhood of
$C_{\varepsilon/2}$.
\begin{lemma}
\label{support}
There exists $\sigma>0$ such that
%
%
\begin{equation}
\label{081901}
\inf_{\nu\in{\mathcal M}_1(B(z,\sigma))}
\liminf_{T\uparrow\infty} Q^T\nu({\tilde
C})>1-\frac{3\varepsilon}{4}.
\end{equation}
In addition, if $\sigma$ is as above, then for any $k\ge1$ and
$t_1,\ldots,t_k\ge0$, we can choose $T_*$ such that
%
%
\begin{equation}\label{081902}
\inf_{\nu\in{\mathcal M}_1(B(z,\sigma))}
Q^{T,t_1,\ldots,t_k}\nu({\tilde C})>1-\frac{3\varepsilon}{4}\qquad
\forall T\ge T_*.
\end{equation}
\end{lemma}
\begin{pf}
The claim made in (\ref{081901}) follows if we can show that there
exists $\sigma>0$ such that
%
%
\begin{equation}\label{Th4}
\liminf_{T\to+\infty} Q^T(y,{\tilde
C})>1-\frac{3\varepsilon}{4} \qquad\forall y\in B (z,\sigma).
\end{equation}
%
To prove (\ref{Th4}), suppose that $\psi$ is a Lipschitz function
such that $\mathbf{1}_{C_{\varepsilon/2}}\le\psi\le\mathbf
{1}_{\tilde C}$.
Since $(P_t\psi)_{t\ge0}$ is equicontinuous at $z$, we can find
$\sigma>0$ such that $|P_t\psi(y)-P_t\psi(z)|<\varepsilon/4$ for
all $y\in B
(z,\sigma)$. We then have
\[
Q^T(y,{\tilde C})\ge\int_\mathcal X\psi(y')Q^T(y,\dd y')\ge\int
_\mathcal X
\psi(y')Q^T(z,\dd y')-\frac{\varepsilon}{4}
\]
and, using (\ref{Th3}), we conclude that
%
%
\begin{equation}\label{Th4b}
\liminf_{T\uparrow\infty} Q^T(y,{\tilde C}) \ge
\liminf_{T\uparrow\infty} Q^T(z,C_{\varepsilon/2})-\frac
{\varepsilon}{4}>
1-\frac{3\varepsilon}{4}.
\end{equation}
Estimate (\ref{081902}) follows directly from (\ref{081901}) and
Lemma \ref{Lemma}.
\end{pf}

Let us return to the proof of Theorem \ref{Theorem}. Let $\sigma>0$
be as in the above lemma and let $\gamma>0$ denote the supremum of
all sums $\alpha_1+\cdots+\alpha_k$ such that there exist
$\nu_1,\ldots, \nu_k\in\mathcal M_1(B (z,\sigma))$ and
%
%
\begin{equation}
Q^{t_1^0,\ldots, t_{m_0}^0}(x)\ge\alpha_1 Q^{t_1^1,\ldots,
t_{m_1}^1}\nu_1+\cdots+ \alpha_k Q^{t_1^k, \ldots, t_{m_k}^k}\nu_k
\end{equation}
for some $t_1^0,\ldots, t_{m_0}^0, \ldots, t_1^k,\ldots, t_{m_k}^k>0$.
In light of Lemma \ref{support}, to deduce (\ref{Th3b}), it is enough
to show that $\gamma>1-\varepsilon/4$.
Assume, therefore, that
%
%
\begin{equation}\label{081903}
\gamma\le1-\frac{\varepsilon}4.
\end{equation}
Let $D$ be a bounded
subset of $X$, let $T_*>0$ be such that
%
%
\begin{equation}\label{Th888}
Q^T(x,D)>1-\frac{\varepsilon}{8} \qquad\forall T\ge T_*,
\end{equation}
and let
%
%
\begin{equation}\label{Th888b}
\alpha:= \inf_{x\in D} \liminf_{T\uparrow\infty}Q^T(x,B
(z,\sigma))>0.
\end{equation}
Let $\alpha_1,\ldots, \alpha_k>0$, $t_1^0,\ldots, t_{m_0}^0,
\ldots,t^k_1,\ldots, t_{m_k}^k>0$ and $\nu_1,\ldots,
\nu_k\in\mathcal M_1(B(z, \sigma))$ be such that
\[
Q^{t_1^0,\ldots, t_{m_0}^0}(x)\ge\alpha_1 Q^{t_1^1,\ldots,
t_{m_1}^1}\nu_1+\cdots+ \alpha_k Q^{t_1^k, \ldots, t_{m_k}^k}\nu_k
\]
and
%
%
\begin{equation}\label{081904}
\gamma-(\alpha_1+\cdots+\alpha_k)<\frac{\alpha\varepsilon
}{64}.
\end{equation}
For a given $t\ge0$, we let
\[
\mu_t:=Q^{t, t_1^0,\ldots, t_{m_0}^0}(x)-\alpha_1 Q^{t,
t_1^1,\ldots, t_{m_1}^1}\nu_1-\cdots- \alpha_k Q^{t, t_1^k, \ldots,
t_{m_k}^k}\nu_k.
\]
By virtue of Lemma \ref{Lemma}, we can choose $T_*>0$ such that $\|Q^{t,
t_1^0,\ldots, t_{m_0}^0}(x)-Q^{t}(x)\|_{\mathrm{TV}}<\varepsilon/16$ for $t\ge T_*$.
Thus, from (\ref{Th888}), we obtain that for such $t$,
\begin{eqnarray*}
\mu_{t}(D)&>&Q^{t}(x,D)-\|Q^{t, t_1^0,\ldots,
t_{m_0}^0}(x)-Q^{t}(x)\|_{\mathrm{TV}}-(\alpha_1+\cdots+\alpha_k)
\\
&\ge&
1-\frac{\varepsilon}{8}-\frac{\varepsilon}{16}-\gamma\stackrel
{(\fontsize{8.36}{10}\selectfont{\mbox{\ref{081903}}})}{\ge}
\frac{\varepsilon}{16}.
\end{eqnarray*}
%
However, this means that for $t\ge T_*$,
\begin{eqnarray*}
\mathop{\lim\inf}_{T\uparrow\infty}Q^T\mu_{t} (B
(z,\sigma))&\stackrel{\mathrm{Fatou}\ \mathrm{lem.}}{\ge}&\int_\mathcal X
\mathop{\lim\inf}_{T\uparrow\infty}Q^T (y,B
(z,\sigma))\mu_{t}(\dd y)
\\
&\ge&\int_D \mathop{\lim\inf}_{T\uparrow
\infty}Q^T (y,B (z,\sigma))\mu_{t}(\dd y) \stackrel{
(\fontsize{8.36}{10}\selectfont{\mbox{\ref{Th888b}}})}{\ge} \frac{\alpha\varepsilon}{16}.
\end{eqnarray*}
Choose $T_*>0$ such that
%
%
\begin{equation}\label{Th999}
Q^T\mu_{t} (B (z,\sigma))>\frac{\alpha\varepsilon}{32}
\qquad\forall
t,T\ge T_*.
\end{equation}
Let $\nu(\cdot):= (Q^T\mu_{t}) (\cdot|B(z,\sigma))$. Of course,
$\nu\in{\mathcal M}_1(B(z,\sigma))$. From (\ref{Th999}) and the definitions
of $\nu$, $\mu_t$, we obtain, however, that for $t,T$ as above,
\[
Q^{T, t, t_1^0,\ldots, t_{m_0}^0}(x)\ge\alpha_1 Q^{T, t,
t_1^1,\ldots, t_{m_1}^1}\nu_1+\cdots+ \alpha_k Q^{T, t, t_1^k,
\ldots, t_{m_k}^k}\nu_k+\frac{\alpha\varepsilon}{32}\nu.
\]
Hence, $\gamma\ge\alpha_1+\cdots+\alpha_k+\alpha\varepsilon/32$,
which clearly contradicts (\ref{081904}).

\subsection*{Proof of the weak law of large numbers}
Recall that ${\mathbb{P}}_{\mu}$ is the path measure
corresponding to $\mu$, the initial distribution of $(Z(t))_{t\ge0}$.
Let $\mathbb{E}_{\mu}$ be the corresponding expectation and
$d_*:=\int\psi\,\dd\mu_*$.
It then suffices to show that
%
%
\begin{equation}
\label{090602}
\lim_{T\to+\infty}
\mathbb{E}_{\mu} \biggl[\frac{1}{T}\int_0^T\psi(Z(t))\,\dd t \biggr]=d_*
\end{equation}
and
%
%
\begin{equation}
\label{090603}
\lim_{T\to+\infty}
\mathbb{E}_{\mu} \biggl[\frac{1}{T}\int_0^T\psi(Z(t))\,\dd t \biggr]^2=
d_*^2.
\end{equation}
Equality (\ref{090602}) is an obvious consequence of weak-$^*$ mean ergodicity.
To show (\ref{090603}), observe that
the expression under the limit equals
%
%
\begin{eqnarray}
\label{090604}
&&\frac{2}{T^2}\int_0^T\int_0^t \biggl(\int_{\mathcal X} P^s(\psi
P_{t-s}\psi)\,\dd\mu\biggr)\,\dd t\,\dd s\nonumber\\[-8pt]\\[-8pt]
&&\qquad=\frac{2}{T^2}\int_0^T(T-s)
\biggl(\int_{\mathcal X} P_s(\psi\Psi_{T-s})\,\dd\mu\biggr)\,\dd s,\nonumber
\end{eqnarray}
where
%
%
\begin{equation}
\label{090605}
\Psi_t(x):=\int_{\mathcal X} \psi(y)Q^t(x,\dd y)=\frac{1}{t}\int_0^t
P_s\psi(x)\,\dd s.
\end{equation}
The following lemma then holds.
\begin{lemma}
\label{lm090601}
For any $\varepsilon>0$ and a compact set $K\subset X$, there exists $t_0>0$
such that
%
%
\begin{equation}
\label{090606}
\forall t\ge t_0\qquad \sup_{x\in K} \biggl|\Psi_t(x)-\int
_{\mathcal X
} \psi\,\dd\mu_* \biggr|<\varepsilon.
\end{equation}
\end{lemma}
\begin{pf}
It suffices to show
equicontinuity of $(\Psi_t)_{t\ge0}$ on any compact set $K$.
The proof
then follows from pointwise convergence of $\Psi_t$ to $d_*$ as $t\to
\infty$
and the Arzela--Ascoli theorem. The equicontinuity of the above family
of functions
is a direct consequence of the e-property and a simple covering argument.
\end{pf}

Now, suppose that $\varepsilon>0$. One can find a compact set $K$ such that
%
%
\begin{equation}
\label{090607}
\forall t\ge0\qquad Q^t\mu(K^c)<\varepsilon.
\end{equation}
Then
\begin{eqnarray*}
&& \biggl|\frac{2}{T^2}\int_0^T(T-s) \biggl(\int_{\mathcal X}
P_s(\psi\Psi
_{T-s})\,\dd\mu\biggr)\,\dd s
-\frac{2d_*}{T^2}\int_0^T(T-s) \biggl(\int_{\mathcal X} P_s\psi\,\dd
\mu
\biggr)\,\dd
s \biggr|\\
&&\qquad
\le I+\mathit{I I},
\nonumber
\end{eqnarray*}
where
\[
I:=\frac{2}{T^2}\int_0^T(T-s) \biggl(\int_{\mathcal X} P_s\bigl(\psi(\Psi
_{T-s}-d_*)\mathbf{1}_K\bigr)\,\dd\mu\biggr)\,\dd s
\]
and
\[
\mathit{I I}:=\frac{2}{T^2}\int_0^T(T-s) \biggl(\int_{\mathcal X} P_s\bigl(\psi
(\Psi
_{T-s}-d_*)\mathbf{1}_{K^c}\bigr)\,\dd\mu\biggr)\,\dd s.
\]
According to Lemma \ref{lm090601}, we can find $t_0$ such that (\ref{090606})
holds with the compact set $K$ and $\varepsilon\|\psi\|_{\infty
}^{-1}$. We
then obtain
$
|I|\le\varepsilon.
$
Also, note that
\[
|\mathit{I I}|\le2\|\psi\|_\infty(\|\psi\|_{\infty}+|d_*|)Q^T\mu
(K^c)\stackrel{(\fontsize{8.36}{10}\selectfont{\mbox{\ref{090607}}})}{<}2\varepsilon\|\psi\|_\infty(\|
\psi\|
_{\infty}+|d_*|).
\]
The limit on the right-hand side of (\ref{090603})
therefore
equals
\begin{eqnarray*}
&&\lim_{T\to+\infty}
\frac{2d_*}{T^2}\int_0^T(T-s) \biggl(\int_{\mathcal X} P_s\psi\,\dd\mu
\biggr)\,\dd s
\\
&&\qquad=
\lim_{T\to+\infty}
\frac{2d_*}{T^2}\int_0^T\dd s\int_0^s \biggl(\int_{\mathcal X}
Q^{s'}\psi
\,\dd\mu\biggr)\,\dd s'
=d_*^2.
\end{eqnarray*}

\section[Proof of Theorem 3]{Proof of Theorem \protect\ref{TGeneral}}
\label{SPDEs}

In what follows, we are going to verify the assumptions of Theorem
\ref{Theorem}. First, observe that (\ref{Th2}) follows from (ii) and
Chebyshev's inequality. The e-property implies equicontinuity of
$(P_t\psi, t\ge0)$ at any point for any bounded, Lipschitz function
$\psi$. What remains to be shown, therefore, is condition (\ref{Th1}).
The rest of the proof is devoted to that objective. It will be given
in five steps.
\renewcommand{\theStep}{\Roman{Step}}
\begin{Step}\label{StepI}
We show that we can find a bounded Borel set $B$ and a
positive constant $r^*$ such that
%
%
\begin{equation}\label{Tw1}
\mathop{\lim\inf}_{T\uparrow
\infty}Q^T(x,B)>\frac{1}{2}\qquad \forall x\in\mathcal K+ r^* B(0,1).
\end{equation}
To prove this, observe, by (ii) and Chebyshev's inequality, that for
every $y\in\mathcal K$, there exists a bounded Borel set $B_{y}^0$
such that
$ \liminf_{T\uparrow\infty}Q^T(y,B_{y}^0)>3/4$. Let $B_{y}$ be a
bounded, open set such that $B_{y}\supset B_{y}^0$ and let $\psi\in
C_b(\mathcal X)$ be such that $\mathbf{1}_{B_{y}}\ge\psi\ge\mathbf
{1}_{B_{y}^0}$.
Since $(P_t\psi)_{ t\ge0}$ is equicontinuous at $y$, we can find
$r_{y}>0$ such that $|P_t\psi(x)-P_t\psi(y)|<1/4$ for all $x\in
B(y,r_{y})$ and $t\ge0$. Therefore, we have
\begin{eqnarray*}
\mathop{\lim\inf}_{T\uparrow
\infty}Q^T(x,B_{y})
&\ge&
\mathop{\lim\inf}_{T\uparrow\infty}
\frac{1}{T}\int_0^T {P}_s\psi(x)\,\dd s\\
&\ge&
\mathop{\lim\inf}_{T\uparrow\infty}
\frac{1}{T}\int_0^T {P}_s\psi(y)\,\dd s-\frac{1}{4}\\
&\ge&
\mathop{\lim\inf}_{T\uparrow\infty} Q^T(y,
B^0_y)-\frac{1}{4}
>\frac{1}{2}.
\end{eqnarray*}
Since the attractor is compact, we can find a finite covering
$B(y_i,r_{y_i})$, $i=1,\ldots,N$, of $\mathcal K$. The claim made in
(\ref{Tw1}) therefore holds for $B:=\bigcup_{i=1}^N B_{y_i}$ and
$r^*>0$ sufficiently small so that $\mathcal K+ r^*B(0,1)\subset
\bigcup_{i=1}^NB(y_i,r_{y_i})$.
\end{Step}
\begin{Step}\label{StepII}
Let $B\subset\mathcal X$ be as in Step \ref{StepI}. We prove that
for every bounded Borel set $D\subset\mathcal X$, there exists a
$\gamma>0$
such that
%
%
\begin{equation}\label{Tw3}
\mathop{\lim\inf}_{T\uparrow
\infty}Q^T(x,B)>\gamma\qquad\forall x\in D.
\end{equation}
From the fact that $\mathcal K$ is a global attractor for (\ref{082102}),
for any $r>0$ and a bounded Borel set $D$, there exists an $L>0$ such
that $Y^x(L)\in\mathcal K+\frac{r}{2}B(0,1)$ for all $x\in D$. By
(\ref{ET240}), we have
\[
p(r,D):= \inf_{x\in D}\mathbb{P} \bigl(\|Z^x(L)-Y^x (L)\|
_{\mathcal X}<r/2
\bigr)>0.
\]
We therefore obtain that
%
%
\begin{equation}\label{082201}
P_L\mathbf{1}_{\mathcal K+rB(0,1)}(x)\ge p(r,D) \qquad\forall x\in
D.
\end{equation}
Let $r^*>0$ be the constant given in Step \ref{StepI}. Then
%
%
\begin{eqnarray}\label{Tw7}
&&\mathop{\lim\inf}_{T\uparrow\infty}Q^T(x,B)\nonumber\\[-3pt]
&&\qquad=\mathop{\lim\inf}_{T\uparrow\infty}
\frac{1}{T}\int_0^T {P}_{s+L}\mathbf{1}_{B}(x)\,\dd s\nonumber\\
&&\qquad=\mathop{\lim\inf}_{T\uparrow\infty}
\frac{1}{T}\int_0^T {P}_{s+L}^*\delta_{x}(B)\,\dd s\nonumber\\
&&\qquad=
\mathop{\lim\inf}_{T\uparrow\infty}
\frac{1}{T}\int_0^T \int_{\mathcal X}P_s\mathbf
{1}_{B}(z){P}_{L}^*\delta
_{z}(\dd
z)\,\dd s\nonumber\\
&&\qquad\ge\mathop{\lim\inf}_{T\uparrow\infty}
\frac{1}{T}\int_0^T \int_{\mathcal K+r^*B(0,1)}
P_s\mathbf{1}_{B}(z){P}_{L}^*\delta_{x}(\dd z )\,\dd s\\
&&\hspace*{-12.62pt}\qquad\mathop{\mathop{\ge}\limits^{\mathrm{and}\ \mathrm{Fatou}}}\limits^{\mathrm{Fubini}} \int_{\mathcal K+r^*
B(0,1)}\liminf_{T\uparrow\infty}Q^T(z,B)
{P}_{L}^*\delta_{x}(\dd z )\nonumber\\
&&\hspace*{-4.16pt}\qquad\stackrel{(\fontsize{8.36}{10}\selectfont{\mbox{\ref{Tw1}}})}{\ge}\frac{1}{2}
\int_{\mathcal X}\mathbf{1}_{\mathcal K+r^*B(0,1)}(z){P}_{L}^*\delta
_{x}(\dd
z)\nonumber\\
&&\qquad=\frac{1}{2} {P}_{L}\mathbf{1}_{\mathcal K+r^*
B(0,1)}(x)\nonumber\\
&&\hspace*{-4.16pt}\qquad\stackrel{(\fontsize{8.36}{10}\selectfont{\mbox{\ref{082201}}})}{\ge}
\gamma:=\frac{p(r^*,D)}{2} \qquad\forall x\in D.\nonumber
\end{eqnarray}
\end{Step}
\begin{Step}\label{StepIII}
We show here that for every bounded Borel set
$D\subset\mathcal X$ and any radius $r>0$, there exists a $w>0$ such that
%
%
\begin{equation}\label{Tw8}
\inf_{x\in D}\mathop{\lim\inf}_{T\uparrow
\infty} Q^T\bigl(x,\mathcal K+r B(0,1)\bigr)>w.
\end{equation}
We therefore fix $D\subset\mathcal X$ and $r>0$. From Step \ref{StepII}, we know that
there exist a bounded set $B\subset\mathcal X$ and a positive constant
$\gamma>0$ such that (\ref{Tw3}) holds. By (\ref{ET240}), we have, as
in (\ref{Tw7}),
%
%
\begin{eqnarray}\label{Tw7b}
&&\mathop{\lim\inf}_{T\uparrow\infty}Q^T\bigl(x,\mathcal K+ r
B(0,1)\bigr)\nonumber\\
&&\hspace*{7pt}\qquad=
\mathop{\lim\inf}_{T\uparrow\infty}\frac{1}{T}\int_0^T
\int_{\mathcal X}P_L\mathbf{1}_{\mathcal K+
rB(0,1)}(z){P}_{s}^*\delta_{x}(\dd z)\,\dd s\\
&&\qquad\stackrel{\mathrm{Fubini}}{\ge} \liminf_{T\uparrow
\infty} \int_{B}P_L\mathbf{1}_{\mathcal K+ rB(0,1)}(z)Q^T(x,\dd z).\nonumber
\end{eqnarray}
Using (\ref{082201}), we can further estimate the last right-hand side of
(\ref{Tw7b}) from below by
%
%
\begin{equation}\label{Tw7c}
p(r,D)\liminf_{T\uparrow\infty}
Q^T(x,B)\stackrel{(\fontsize{8.36}{10}\selectfont{\mbox{\ref{Tw3}}})}{>} p(r,D)\gamma.
\end{equation}
We therefore obtain (\ref{Tw8}) with $w=\gamma p(r,D)$.
\end{Step}
\begin{Step}\label{StepIV}
Choose $z\in\bigcap_{y\in\mathcal K}\bigcup_{t\ge0}
\Gamma^t(y)\neq\varnothing$. We are going to show that for every
$\delta>0$, there exist a finite set of positive numbers $S$ and a
positive constant ${\tilde r}$ satisfying
%
%
\begin{equation}\label{Tw2}
\inf_{x\in\mathcal K+\tilde r B(0,1)}\max_{s\in S} P_s\mathbf
{1}_{B(z,\delta
)} (x)>0.
\end{equation}

Let $t_{x}>0$ for $x\in\mathcal K$ be such that $z\in\operatorname{supp}
P_{t_{x}}^*\delta_{x}$. By the Feller property of $(P_t)_{t\ge0}$, we
may find, for any $x\in\mathcal K$, a positive constant $r_{x}$ such that
%
%
\begin{equation}\label{Twier1}
P_{t_{x}}^*\delta_{y} (B(z, \delta))\ge P_{t_{x}}^*\delta_{x} (B(z,
\delta))/2 \qquad\mbox{for $y\in B(x, r_{x})$}.
\end{equation}
Since $\mathcal K$ is compact, we may choose $x_1, \ldots, x_p\in
\mathcal K$ such
that $\mathcal K\subset\bigcup_{i=1}^p B_i$, where $B_i=B(x_i, r_{x_i})$
for $i=1, \ldots, p$.
Choose ${\tilde r}>0$ such that $\mathcal K+{\tilde r} B(0, 1)\subset
\bigcup
_{i=1}^p B_i$.
\end{Step}
\begin{Step}\label{StepV}
Fix a bounded Borel subset $D\subset\mathcal X$,
$z\in\bigcap_{y\in\mathcal K}\bigcup_{t\ge0}\Gamma^t (y)$ and
$\delta
>0$. Let a positive constant $\tilde r$
and a finite set $S$ be such that (\ref{Tw2}) holds. Set
%
%
\begin{equation}\label{082401}
u:=\inf_{x\in\mathcal K+\widehat r B(0,1))}\max_{s\in S}
P_s\mathbf{1}_{B(z,\delta)} (x)>0.
\end{equation}
From Step \ref{StepIII}, it follows that there exists $w>0$ such that
(\ref{Tw8}) holds for $r=\tilde r$.

Denote by $\# S$ the cardinality of $S$. We can easily check that
%
%
\begin{eqnarray}\label{Tw4}
&&\mathop{\lim\inf}_{T\uparrow\infty} \sum_{q\in
S}\frac{1}{T}\int_0^T {P}_{q+s}\mathbf{1}_{B(z,\delta)}(x)\,\dd s \nonumber\\[-8pt]\\[-8pt]
&&\qquad=\#
S\mathop{\lim\inf}_{T\uparrow
\infty}Q^T(x,B(z,\delta)) \qquad\forall x\in D.\nonumber
\end{eqnarray}
On the other hand, we have
%
%
\begin{eqnarray}\label{Tw5}
&&\hspace*{6.5pt}\sum_{q\in S}\frac{1}{T}\int_0^T {P}_{q+s}\mathbf{1}_{B(z,\delta
)}(x)\,\dd s
\nonumber\\
&&\hspace*{6.5pt}\qquad=\int_{\mathcal X}\sum_{q\in S}P_{q}\mathbf{1}_{B(z,\delta
)}(y)Q^T(x,\dd
y)\nonumber\\[-8pt]\\[-8pt]
&&\hspace*{6.5pt}\qquad\ge \int_{\mathcal K+ \tilde rB(0,1)}\sum_{q\in S}P_{q}\mathbf{1}
_{B(z,\delta)}(y)Q^T(x,\dd y)\nonumber\\
&&\qquad\stackrel{(\fontsize{8.36}{10}\selectfont{\mbox{\ref{082401}}})}{\ge}
uQ^T\bigl(x,\mathcal K+ \tilde
rB(0,1)\bigr) \qquad\forall x\in D.\nonumber
\end{eqnarray}
Combining (\ref{Tw8}) with (\ref{Tw5}), we obtain
\[
\mathop{\lim\inf}_{T\uparrow\infty} \sum_{q\in
S}\frac{1}{T}\int_0^T {P}_{q+s}\mathbf{1}_{B(z,\delta)}(x)\,\dd s
>uw \qquad\forall x\in D,
\]
and, finally, by (\ref{Tw4}),
\[
\mathop{\lim\inf}_{T\uparrow
\infty}Q^T(x,B(z,\delta))>uw/\#S\qquad \forall x\in D.
\]
This shows that condition (\ref{Th1}) is satisfied with
$\alpha=uw/\#S$.
\end{Step}

\section{Ergodicity of the Lagrangian observation process}\label
{Observable-section}

This section is in preparation for the proof of Theorem \ref{TTracer}.
Given an $r\ge0$, we denote by $\mathcal X^r$ the Sobolev space which
is the
completion of
\[
\biggl\{ x \in C^\infty(\mathbb{T}^d;\mathbb{R}^d)\dvtx\int
_{\mathbb{T}^d}
x(\xi)\,\dd\xi=0, \widehat x(k)\in\operatorname{Im} \mathcal E(k),
\forall k\in\mathbb{Z}^d_* \biggr\}
\]
with respect to the norm
\[
\|x\|^2_{\mathcal X^r}: =\sum_{k\in\mathbb{Z}^d_*} |k|^{2r}|\widehat x(k)|^2,
\]
where
\[
\widehat x(k):=(2\pi)^{-d}\int_{\mathbb T^d}x(\xi)e^{-i\xi\cdot
k}\,\dd
\xi,\qquad k\in\mathbb{Z}^d,
\]
are the Fourier coefficients of $x$. Note that
$\mathcal X^u\subset\mathcal X^r$ if $u>r$.

Let $A_r$ be an operator on $\mathcal X^r$ defined by
%
%
\begin{equation}
\label{A-r} \widehat{A_r x}(k):= -\gamma(k)\widehat x(k),\qquad
k\in\mathbb{Z}^d_*,
\end{equation}
with the domain
%
%
\begin{equation}
\label{D-r} D(A_r):= \biggl\{x \in\mathcal X^r\dvtx\sum_{k\in
\mathbb{Z}^d_*}
|\gamma(k)|^2 |k|^{2r} \vert\widehat x (k) \vert^2
<\infty\biggr\}.
\end{equation}
Since the operator is self-adjoint, it generates a $C_0$-semigroup
$(S_r(t))_{t\ge0}$ on $\mathcal X^r$. Moreover, for $u>r$, $A_u$ is
the restriction
of $A_r$ and $S_u$ is the restriction of $S_r$. From now on, we will
omit the subscript $r$ when it causes no confusion, writing $A$ and
$S$ instead of $A_r$ and $S_r$, respectively.

Let $Q$ be a symmetric positive definite bounded linear operator on
\[
\biggl\{x\in L^2(\mathbb{T}^d,\dd\xi;\mathbb{R}^d)\dvtx
\int_{\mathbb{T}^d}x(\xi)\,\dd\xi=0 \biggr\}
\]
given by
\[
\widehat{Qx}(k) := \gamma(k) \mathcal E(k)\widehat x(k), \qquad
k\in \mathbb{Z}^d_*.
\]
Let $m$ be the constant appearing in (\ref{H1}) and let
$\mathcal X:=\mathcal X^m$ and $\mathcal V:= \mathcal X^{m+1}$. Note
that, by Sobolev embedding
(see, e.g., Theorem 7.10, page 155 of \cite{Trudinger}),
$\mathcal X\hookrightarrow C^1(\mathbb{T}^d,\mathbb{R}^d)$ and hence
there exists a
constant $C>0$ such that
%
%
\begin{equation}\label{EA1}
\|x\|_{C^1(\mathbb{T}^d;\mathbb{R}^d)}\le C\|x\|_{\mathcal X}
\qquad\forall x\in{\mathcal X}.
\end{equation}

For any $t>0$, the operator $S(t)$ is bounded from any $\mathcal X^r$ to
$\mathcal X^{r+1}$. Its norm can be easily estimated by
\[
\|S(t)\|_{L(\mathcal X^r,\mathcal X^{r+1})}\le\sup_{k\in\mathbb
{Z}^d_*}|k|e^{-\gamma(k)t}.
\]
Let $e_k(x):= e^{ik\cdot x}$, $k\in\mathbb{Z}^d$. The
Hilbert--Schmidt norm
of the operator $S(t)Q^{1/2}$ (see Appendix C of \cite{DZ1}) is
given by
\begin{eqnarray*}
\| S(t)Q^{1/2}\| ^2_{L_{(\mathrm{HS})}(\mathcal X,\mathcal V)}:\!&=& \sum
_{k\in\mathbb{Z}^d}\|
S(t)Q^{1/2}e_k\| ^2_{\mathcal V}\\
&=& \sum_{k\in\mathbb{Z}^d}|k|^{2(m+1)}
\gamma(k)e^{-2\gamma(k)t}\operatorname{Tr} {\mathcal E}(k).
\end{eqnarray*}
Taking into account assumptions (\ref{H1}) and (\ref{H2}),
we easily obtain the following lemma.
\begin{lemma}\label{L1}
\textup{(i)} For each $t>0$, the operator $Q^{1/2}S(t)$ is Hilbert--Schmidt
from $\mathcal X$ to
$\mathcal V$ and there exists $\beta\in(0,1)$ such that
\[
\int_0^\infty t^{-\beta}\| S(t)Q^{1/2}\| ^2_{L_{(\mathrm{HS})}(\mathcal
X,\mathcal V)} \,\dd t <\infty.
\]

\textup{(ii)} For any $r\ge0$ and $t>0$, the operator $S(t)$ is bounded
from $\mathcal X^r$ into $\mathcal X^{r+1}$ and
\[
\int_0^\infty\|S(t)\|_{L(\mathcal X^r,\mathcal X^{r+1})} \,\dd t<\infty.
\]
\end{lemma}

Let $W=(W(t))_{ t\ge0}$ be a cylindrical Wiener process in $\mathcal X$
defined on a filtered probability space $\mathfrak{A}=
(\Omega,\mathcal F,(\mathcal F_t),\mathbb{P})$. By Lemma \ref{L1}(i)
and Theorem 5.9,
page 127 of \cite{DZ1}, for any $x\in\mathcal X$, there exists a
unique, continuous
in $t$, $\mathcal X$-valued process $V^x$ solving, in the mild sense, the
Ornstein--Uhlenbeck equation
%
%
\begin{equation}\label{EOU}
\dd V^x(t)= A V^x(t)\,\dd t+ Q^{1/2} \,\dd W(t),\qquad V^x(0)=x.
\end{equation}
Moreover, (\ref{EOU}) defines a Markov family on $\mathcal X$ (see Section
9.2 of \cite{DZ1}) and the law $\mathcal L(V(0,\cdot))$ of $V(0,\cdot
)$ on
$\mathcal X$ is its unique invariant probability measure (see Theorem 11.7
of~\cite{DZ1}). Note that, since $m>d/2+1$, for any fixed $t$, the
realization of $V^x(t,\xi)$ is Lipschitz in the $\xi$ variable. If
the filtered probability space $\mathfrak{A}$ is sufficiently
rich, that is, if there exists an $\mathcal F_0$-measurable random
variable with law $\mathcal L(V(0,\cdot))$, then the stationary solution
to (\ref{EOU}) can be found as a stochastic process over
$\mathfrak{A}$. Its law on the space of trajectories
$C([0,\infty)\times\mathbb{T}^d;\mathbb{R}^d)$ coincides with the
law of
$(V(t,\cdot))_{ t\ge0}$.

\subsection{An evolution equation describing the environment
process}
\label{sec6}

Since the realizations of $V^x(t,\cdot)$ are Lipschitz in the spatial
variable, equation (\ref{ET1}), with $V^x(t,\xi)$ in place of
$V(t,\xi)$, has a unique solution $\mathbf x_x(t)$, $t\ge0$, for given
initial data $\mathbf x_0$. In fact, with no loss of generality, we
may, and
shall, assume that $\mathbf x_0=0$. In what follows, we will also
denote by
$\mathbf x$ the solution of (\ref{ET1}) corresponding to the stationary
right-hand side $V$. Let $\mathcal Z(s,\xi):=V(s,\xi+\mathbf x(s))$ be
the \textit{Lagrangian observation of the environment process} or, in
short, \textit{the observation process}. It is known (see \cite{FKP}
and \cite{KP}) that $\mathcal Z(s,\cdot)$ solves the equations
%
%
\begin{eqnarray}\label{Zet}
\dd\mathcal Z(t) &=& [A\mathcal Z(t) + B(\mathcal Z(t),\mathcal
Z(t)) ]\,\dd t +
Q^{1/2}\,\dd
\tilde W(t),\nonumber\\[-8pt]\\[-8pt]
\mathcal Z(0,\cdot) &=& V\bigl(0, \mathbf x(0)+\cdot\bigr),\nonumber
\end{eqnarray}
where $\tilde W$ is a certain cylindrical Wiener process on the
original probability space $\mathfrak{A}$ and
%
%
\begin{eqnarray}\label{form-B}
B(\psi,\phi)(\xi):= \Biggl( \sum_{j=1}^d \psi_j(0)\,\frac{\partial
\phi_1}{\partial\xi_j}(\xi),\ldots, \sum_{j=1}^d
\psi_j(0)\,\frac{\partial\phi_d}{\partial\xi_j}(\xi)
\Biggr),\nonumber\\[-8pt]\\[-8pt]
\eqntext{\psi,\phi\in\mathcal X, \xi\in\mathbb{T}^d.}
\end{eqnarray}
By (\ref{EA1}), $B(\cdot,\cdot)$ is a continuous bilinear form mapping
from $\mathcal X\times\mathcal X$ into $\mathcal X^{m-1}$.

For a given an ${\mathcal F}_0$-measurable random variable $Z_0$ which
is square-integrable in
$\mathcal X$ and a cylindrical Wiener process $W$ in
$\mathcal X$, consider the SPDE
%
%
\begin{equation}\label{EW2}\qquad
\dd Z(t) = [AZ(t) + B(Z(t),Z(t)) ]\,\dd t + {Q}^{1/2} \,\dd
W(t),\qquad Z(0)=Z_0.
\end{equation}
Taking into account Lemma \ref{L1}(ii), the local existence and
uniqueness of a mild solution follow by a standard Banach fixed
point argument. For a different type of argument, based on the Euler
approximation scheme, see Section 4.2 of \cite{FKP}. Global
existence also follows; see the proof of the moment estimates in
Section~\ref{sec5.3.2} below.

Given $x\in\mathcal X$, let $Z^x(t)$ denote the value at $t\ge0$ of a
solution to (\ref{EW2}) satisfying $Z^x(0,\xi) =x(\xi), \xi\in
\mathbb{T}^d$. Since the existence of a solution follows from the Banach
fixed point argument, $Z=(Z^x, x\in\mathcal X)$ is a stochastically
continuous Markov family and its transition semigroup $(P_t)_{t\ge0}$ is
Feller; for details see, for example, \cite{DZ1} or \cite{PZ}. Note that
\[
P_t\psi(x):=\mathbb{E} \psi\bigl(V^x\bigl(t,\mathbf x_x(t)+\cdot\bigr)\bigr).
\]

The following result on ergodicity of the observation process,
besides being of independent interest, will be crucial for the proof
of Theorem \ref{TTracer}.
\begin{theorem}\label{Pttheorem} Under assumptions (\ref{H1}) and
(\ref{H2}), the
transition semigroup $(P_t)_{t\ge0}$ for the family $Z=(Z^x, x\in
\mathcal X)$ is
weak-$^*$ mean ergodic.
\end{theorem}

To prove the above theorem, we verify the hypotheses of Theorem
\ref{TGeneral}.

\subsubsection{Existence of a global attractor} Note that
$Y^0(t)\equiv0$ is the global attractor for the semi-dynamical system
$Y= (Y^x, x\in\mathcal X )$ defined by the deterministic
problem
%
%
\begin{equation}\label{ET21}
\frac{\dd Y^x(t)}{\dd t}= A Y^x(t) + B(Y^x(t),Y^x(t)),\qquad Y^x(0)=x.
\end{equation}
Clearly, this guarantees the uniqueness of an invariant measure
$\nu_*$ for the corresponding semi-dynamical system; see Definition
\ref{def4.2}. Our claim follows from the exponential stability of
$Y^0$, namely,
%
%
\begin{equation}\label{ET22}
\forall x\in\mathcal X, t>0\qquad \|Y^x(t) \|_{\mathcal X}\le
e^{-\gamma_*
t}\|x\|_{\mathcal X},
\end{equation}
where
%
%
\begin{equation}\label{lgamma}
\gamma_* = \inf_{k\in\mathbb{Z}^d_*}\gamma(k)
\end{equation}
is strictly positive by (\ref{H2}). Indeed, differentiating $\|Y^x(t)
\|_{\mathcal X}^2$ over $t$, we obtain
\[
\frac{\dd}{\dd t}\| Y^x (t) \|_{\mathcal X}^2=2\langle A Y^x(t),Y^x
(t)\rangle_{\mathcal X}+2\sum_{j=1}^d Y^x(t,0) \biggl\langle\frac
{\partial
Y^x(t)}{\partial\xi_j}, Y^x(t) \biggr\rangle_{\mathcal X}.
\]
The last term on the right-hand side vanishes, while the first one
can be estimated from above by $-2\gamma_*\|Y^x(t) \|_{\mathcal X}^2$.
Combining these observations with Gronwall's inequality, we obtain
(\ref{ET22}).

\subsubsection{Moment estimates}
\label{sec5.3.2} Let $B(0,R)$ be the ball in $\mathcal X$ with center at
$0$ and radius $R$. We will show that for any $R>0$ and any integer
$n\ge1$,
%
%
\begin{equation}\label{ET23}
\sup_{x\in B(0, R)} \sup_{t\ge0}\mathbb{E} \| Z^x(t)\|_{\mathcal
X}^{2n}<\infty.
\end{equation}
Recall that $V^x$ is the solution to (\ref{EOU}) satisfying
$V^x(0)=x$. Let $\mathbf x_x= (\mathbf x_x(t), t\ge0 )$
solve the
problem
%
%
\begin{equation}\label{Ex1}
\frac{\dd\mathbf x_x }{\dd t}(t) = V^x(t,\mathbf x_x(t)),\qquad \mathbf
x_x(0)=0.
\end{equation}
We then obtain
%
%
\begin{equation}
\label{082901b} \|Z^x(t)\|_{\mathcal X}^2\stackrel{\dd}{=}\int
_{\mathbb{T}^d} \bigl\vert\nabla^m
V^x\bigl(t,\mathbf x_x(t)+\xi\bigr) \bigr\vert^2\,\dd\xi=\|V^x(t)\|_{\mathcal X}^{2},
\end{equation}
where the first equality means equality in law.
Since $V^x$ is Gaussian, there is a constant $C_1>0$ such that
\[
\mathbb{E} \|V^x(t)\|_{\mathcal X}^{2n}\le C_1 (\mathbb{E}
\|V^x(t)\|_{\mathcal X}^2 )^n.
\]
Hence, there is a constant $C_2>0$ such that for $\|x\|_\mathcal X\le R$,
\begin{eqnarray*}
\mathbb{E} \|V^x(t)\|_{\mathcal X}^{2n}&\le& C_2 (1+R^{2n}
) \biggl(
\int_0^t \|
S(t-s)Q^{1/2}\|_{L_{(\mathrm{HS})}(\mathcal X,\mathcal X)}^2\,\dd s \biggr)^{n} \\
&\le& C_{2} (1 +R^{2n} ) \biggl( \int_0^\infty\|
S(s)Q^{1/2}\|_{L_{(\mathrm{HS})}(\mathcal X,\mathcal X)}^2\,\dd s \biggr)^{n}.
\end{eqnarray*}
Note that there is a constant $C_3$ such that
\[
\int_0^\infty\| S(s)Q^{1/2}\|_{L_{(\mathrm{HS})}(\mathcal X,\mathcal X)}^2\,\dd s
\le C_3
\nnorm\mathcal E\nnorm^2<\infty,
\]
where $\nnorm\mathcal E\nnorm^2$ appears in (\ref{H1}) and (\ref{ET23})
indeed follows.

\subsubsection{Stochastic stability}

Define $\tilde Z^x(t):=V^x(t,\mathbf x_x(t)+\cdot)$. This satisfies
equation (\ref{Zet}) and so the laws of $(\tilde Z^x(t))_{t\ge0}$ and
$( Z^x(t))_{t\ge0}$ are identical.
On the other hand, for $\tilde V^x(t):=S(t)x$ and
\[
\frac{\dd\mathbf y_x}{\dd t}=\tilde V^x(t,\mathbf y_x(t)),\qquad \mathbf y(0)=0,
\]
we have that $Y^x(t,\cdot):=\tilde V^x(t,\mathbf y_x(t)+\cdot)$ satisfies
(\ref{ET21}).
To show stochastic stability, it suffices to prove that
%
%
\begin{equation}\label{ET24}
\forall \varepsilon,  R,  T>0\qquad
\inf_{\|x\|_{\mathcal X}\le R}\mathbb{P} \bigl(\|\tilde Z^x(T)-Y^x
(T)\|_{\mathcal X}<\varepsilon\bigr)>0.
\end{equation}
%

Let $M=(M(t))_{ t\ge0}$ be
the stochastic convolution process
%
%
\begin{equation}
\label{082903}
M(t) := \int_0^t S(t-s) Q^{1/2} \,\dd W(s),\qquad t\ge0.
\end{equation}
It is a centered, Gaussian, random element in the Banach space
$C([0,T],\mathcal X)$ whose norm we denote by $\|\cdot\|_\infty$. We will
use the same notation for the norm on $C[0,T]$.
Note that $V^x(t)=\tilde V^x(t)+M(t)$.
Since $M$ is a centered, Gaussian, random element in the Banach
space $C([0,T],\mathcal X)$, its topological support is a
closed linear subspace; see, for example, \cite{VAK}, Theorem 1, page
61. Thus,
in particular, $0$~belongs to the support of its law and
%
%
\begin{equation}\label{ET25}
\forall \delta>0\qquad q:=\mathbb{P} ( F_\delta) >0,
\end{equation}
where $F_\delta:=[\|M\|_{\infty}
<\delta]$.
Since $\|V^x-\tilde V^x\|_{\infty}<\delta$ on $F_\delta$, we can
choose $\delta$ sufficiently small so that
$
\|\mathbf x_x-\mathbf y_x\|_{\infty}<\rho,
$
where $\rho$ is chosen in such a way that
\begin{eqnarray*}
&&\|\tilde Z^x(T)-Y^x(T)\|_{\mathcal X}\\
&&\qquad\le\bigl\|V^x\bigl(T,\mathbf x_x(T)+\cdot
\bigr)-\tilde
V^x\bigl(T,\mathbf x_x(T)+\cdot\bigr)\bigr\|_{\mathcal X}\\
&&\qquad\quad{} +\bigl\|\tilde V^x\bigl(T,\mathbf x_x(T)+\cdot\bigr)-\tilde V^x\bigl(T,\mathbf
y_x(T)+\cdot\bigr)\bigr\|
_{\mathcal X}<\varepsilon\qquad\forall x\in B(0,R)
\end{eqnarray*}
on $F_\delta$.
Hence, (\ref{ET24}) follows.

\subsubsection{e-property of the transition semigroup}

It suffices to show that for any $\psi\in C^1_b(\mathcal X)$ and $R>0$,
there exists a positive constant $C$ such that
%
%
\begin{equation}\label{ET26}
\sup_{t\ge0}\sup_{\|x\|_{\mathcal X}\le R}\|D P_t\psi(x)\|
_{\mathcal X}\le C \|
\psi\|_{C_b^1({\mathcal X})}.
\end{equation}
Here, $D\phi$ denotes the Fr\'echet derivative of a given function
$\phi\in C^1_b(\mathcal X)$.
Indeed, let $\rho_n\in C^2_0(\mathbb{R}^n)$ be supported in the ball of
radius $1/n$, centered at $0$ and such that $\int_{\mathbb{R}^n}\rho
_n(\xi)\,\dd\xi=1$. Suppose that $(e_n)$ is an orthonormal base in $\mathcal
X$ and
$Q_n$ is the orthonormal projection onto $\operatorname{span}\{e_1,\ldots,e_n\}$.
Define
\[
\psi_n(x):=\int_{\mathbb{R}^n}\rho_n(Q_n x-\xi)\psi\Biggl(\sum
_{i=1}^n\xi_ie_i \Biggr)\,\dd\xi,\qquad x\in\mathcal X.
\]
One can deduce (see part 2 of the proof of Theorem 1.2, pages 164 and 165 in
\cite{PZ1}) that for any $\psi\in\operatorname{Lip} (\mathcal X)$, the
sequence $(\psi_n)$ satisfies $(\psi_n)\subset C_b^1({\mathcal X})$ and
$\lim_{n\to\infty}\psi_n(x)=\psi(x)$ pointwise. In addition, $\|
\psi_n\|_{L^\infty}\le\|\psi\|_{L^\infty}$ and ${\sup_{z}}\|D \psi
_n(z)\|_{\mathcal X}\le\operatorname{Lip}(\psi)$. Let $R>0$ be arbitrary
and $x,y\in B(0,R)$. We can write
\begin{eqnarray*}
|P_t\psi(x)-P_t\psi(y)|&=&\lim_{n\to\infty}|P_t\psi_n(x)-P_t\psi
_n(y)|\\
&\le&
\sup_{\|z\|_{\mathcal X}\le R}\|D P_t\psi_n(z)\|_{\mathcal
X}\|x-y\|
_{\mathcal X}\\
&\stackrel{(\fontsize{8.36}{10}\selectfont{\mbox{\ref{ET26}}})}{\le}& C \|\psi_n\|_{C_b^1({\mathcal X})}
\|x-y\|
_{\mathcal X}\\
&\le&
C [\|\psi\|_{\infty}+\operatorname{Lip}(\psi)] \|
x-y\|
_{\mathcal X}.
\end{eqnarray*}
This shows equicontinuity of $(P_t\psi)_{t\ge0}$ for an arbitrary
Lipschitz function $\psi$ in the neighborhood of any $x$ and the
e-property follows.

To prove (\ref{ET26}), we adopt the method from \cite{HM}. First,
note that $
D P_t \psi(x)[v]$, the value of $D P_t\psi(x)$ at $v\in\mathcal X$, is
equal to $ \mathbb{E} \{ D\psi(Z^x(t)) [U(t)] \}$, where
$U(t):=\partial Z^x(t)[v]$ and
\[
\partial Z^x(t)[v] := \lim_{\varepsilon\downarrow0}
\frac{1}{\varepsilon} \bigl( Z^{x+\varepsilon v}(t)- Z^x(t) \bigr),
\]
the limit here taken in $L^2(\Omega,\mathcal F,\mathbb{P};\mathcal
X)$. The process
$U= (U(t), t\ge0 )_{t\ge0}$ satisfies the linear evolution
equation
%
%
\begin{eqnarray}
\label{ET27}
\frac{\dd U(t)}{\dd t}&=& A U(t)+ B(Z^x(t), U(t)) + B(U(t),
Z^x(t)),\nonumber\\[-8pt]\\[-8pt]
U(0)&=&v.\nonumber
\end{eqnarray}
Suppose that $H$ is a certain Hilbert space and $\Phi\dvtx
\mathcal X\to H$ a Borel measurable function. Given an
$(\mathcal F_t)_{t\ge0}$-adapted process $g\dvtx[0,\infty) \times
\Omega
\to
\mathcal X$ satisfying $ \mathbb{E}\int_0^t \|g_s\|^2_{\mathcal X}\,\dd
s <\infty$ for each
$t\ge0$, we denote by $\mathcal D_g\Phi(Z^x(t))$ the Malliavin derivative
of $\Phi(Z^x(t))$ in the direction of $g$. That is, the
$L^2(\Omega,\mathcal F,\mathbb{P};H)$-limit, if exists, of
\[
\mathcal D_g\Phi(Z^x(t)):=\lim_{\varepsilon\downarrow
0}\frac{1}{\varepsilon} [ \Phi(Z^x_{\varepsilon g}(t)) -\Phi(
Z^x(t)) ],
\]
where $Z^x _{ g}(t)$, $t\ge0$, solves the equation
\begin{eqnarray*}
\dd Z^x_g(t) &=& [AZ^x_g(t)+ B(Z^x_g(t),Z^x_g(t)) ] \,\dd t  +
Q^{1/2} \bigl( \dd W(t) + g_t\,\dd t \bigr),\\
Z^x_g(0) &=& x.
\end{eqnarray*}
In particular, one can easily show that when $H=\mathcal X$ and $\Phi=I$,
where $I$ is the identity operator,
the Malliavin derivative of $Z^x(t)$ exists and
the process $D(t):= \mathcal D_gZ^x(t)$, $t\ge0$, solves the linear equation
%
%
\begin{eqnarray}\label{ET28}\quad
\frac{\dd D}{\dd t}(t)&=&A D(t) + B(Z^x(t),D(t))
+B(D(t),Z^x(t)) + Q^{1/2}g(t),\nonumber\\[-8pt]\\[-8pt]
D(0)&=&0.\nonumber
\end{eqnarray}
The following two facts about the Malliavin derivative will be crucial
for us in the sequel. Directly from the definition of the Malliavin
derivative, we derive \textit{the chain rule}: if we suppose that $\Phi
\in
C^1_b(\mathcal X;H)$, then
%
%
\begin{equation}\label{083002}
\mathcal D_g\Phi( Z^x(t))=D\Phi(Z^x(t))[D(t)].
\end{equation}
In addition, the \textit{integration by parts formula} holds; see
Lemma 1.2.1, page 25 of~\cite{Nualart}. If we suppose that $\Phi\in
C^1_b(\mathcal X)$, then
%
%
\begin{equation}
\label{083003} \mathbb{E}[\mathcal D_g\Phi( Z^x(t))]=
\mathbb{E} \biggl[\Phi(Z^x(t))\int_0^t\langle g(s),Q^{1/2}\,\dd
W(s)\rangle_{\mathcal X} \biggr].
\end{equation}
We also have the following proposition.
\begin{proposition}
\label{m-lm} For any given $v,x\in\mathcal X$ such that $\|v\|
_{\mathcal X}\le
1,
\|x\|_{\mathcal X}\le R$, one can find an $(\mathcal F_t)$-adapted
$\mathcal X$-valued
process $g_t= g_t(v,x)$ that satisfies
%
%
\begin{eqnarray}\label{ET29}
\sup_{\|v\|_{\mathcal X}\le1}\sup_{\|x\|_{\mathcal X}\le R} \int
_0^\infty\mathbb{E}
\|Q^{1/2}g_s\|^2_{\mathcal X}\,\dd s &<&\infty,
\\
%
%
\label{ET210}
\sup_{\|v\|_{\mathcal X}\le1}\sup_{\|x\|_{\mathcal X}\le R}\sup
_{t\ge0} \mathbb{E}
\|D Z^x(t)[v]- {\mathcal D}_gZ^x(t)\|_{\mathcal X}&<&\infty.
\end{eqnarray}
\end{proposition}

We prove this proposition shortly. First, however, let us demonstrate
how it can be used to complete the argument for the e-property. Let
$\omega_t(x):={\mathcal D}_gZ^x(t)$ and $ \rho_t(v,x):=D Z^x(t)[v]-
{\mathcal D}_gZ^x(t)$. Then
\begin{eqnarray*}
D P_t\psi(x)[v]&=& \mathbb{E} \{ D\psi(Z^x(t))[\omega
_t(x)] \}
+\mathbb{E} \{ D
\psi(Z^x(t))[\rho_t(v,x)] \}\\
& = &\mathbb{E} \{{\mathcal D}_g \psi(Z^x(t)) \}
+\mathbb{E} \{
D\psi(Z^x(t))[\rho_t(v,x)] \}
\\
&\stackrel{(\fontsize{8.36}{10}\selectfont{\mbox{\ref{083003}}})}{=}& \mathbb{E}
\biggl\{\psi(Z^\xi(t))\int_0^t\langle g(s),Q^{1/2}\,\dd
W(s)\rangle_{\mathcal X} \biggr\}\\
&&{} + \mathbb{E} \{
D\psi(Z^x(t))[\rho_t(v,x)] \}.
\end{eqnarray*}
We have
\[
\biggl| \mathbb{E} \biggl\{\psi(Z^x(t))\int_0^t\langle
g(s),Q^{1/2}\,\dd
W(s)\rangle_{\mathcal X} \biggr\} \biggr|\le\|\psi\|_{L^\infty
} \biggl(\mathbb{E}
\int_0^\infty\|Q^{1/2}g(s)\|^2_{\mathcal X}\,\dd s \biggr)^{1/2}
\]
and
\[
|\mathbb{E} \{
D\psi(Z^x(t))[\rho_t(v,x)] \} |\le\|\psi\|
_{C_b^1({\mathcal X
})}\mathbb{E}
\|\rho_t(v,x)\|_{\mathcal X}.
\]
Hence, by (\ref{ET29}) and (\ref{ET210}), we derive the desired
estimate (\ref{ET26}) with
\begin{eqnarray*}
C &=& \biggl(\mathbb{E}
\int_0^\infty\|Q^{1/2}g(s)\|^2_{\mathcal X}\,\dd s \biggr)^{1/2}\\
&&{}+\sup
_{\|v\|
_{\mathcal X}\le1}\sup_{\|x\|_{\mathcal X}\le R}\sup_{t\ge0}
\mathbb{E}
\|D Z^x(t)[v]- {\mathcal D}_gZ^x(t)\|_{\mathcal X}.
\end{eqnarray*}
Therefore, the e-process property would be shown
if we could prove Proposition~\ref{m-lm}.

\subsubsection[Proof of Proposition 2]{Proof of Proposition \protect\ref{m-lm}} Let us
denote by $\Pi_{\ge N}$
the orthogonal projection onto $\operatorname{span}\{z e^{\ii k\xi}\dvtx
|k|\ge N, z\in\operatorname{Im} \mathcal E(k) \}$ and let $\Pi_{<
N}:=I-\Pi_{\ge N}=\Pi_{\ge N}^\bot$. Write
\[
A_N:=\Pi_{\ge N}A,\qquad Q_N:=\Pi_{\ge N}Q,\qquad A_N^\perp:=\Pi_{<
N}A,\qquad Q_N^\perp:=\Pi_{< N}Q.
\]

Given an integer $N$, let $\zeta^N(v,x)(t)$ be the solution of the problem
%
%
\begin{eqnarray}\label{ET211}\qquad
\frac{\dd\zeta^N}{\dd t}(t)&=&A_N\zeta^N(t)+ \Pi_{\ge N} \bigl(B(Z^x(t),
\zeta^N(t)) +
B(\zeta^N(t), Z^x(t)) \bigr)\nonumber\\
&&{} -\frac12\Pi_{<N}\zeta^N(t)\|\Pi_{<N}\zeta^N(t)\|
_{\mathcal X
}^{-1},\\
\zeta^N(0)&=&v.\nonumber
\end{eqnarray}
We adopt the convention that
%
%
\begin{equation}\label{ET212}
\Pi_{<N}\zeta^N\|\Pi_{<N}\zeta^N\|_{\mathcal X}^{-1}:=0 \qquad
\mbox{if $\Pi_{<N}\zeta^N=0$.}
\end{equation}
Let
%
%
\begin{equation}\label{083005}
g:=Q^{-1/2}f,
\end{equation}
where
%
%
\begin{eqnarray}\label{ET213}\quad
f(t)&:=& A_N^{\perp}\zeta^N(v,x)(t)\nonumber\\
&&{}+ \Pi_{<N} [
B(Z^x(t),\zeta^N(v,x)(t)) + B(\zeta^N(v,x)(t),Z^x(t))
]\\
&&{}
+\tfrac12\Pi_{<N}\zeta^N(v,x)(t)\|\Pi_{<N}\zeta^N(v,x)(t)\|
_{\mathcal X}^{-1}\nonumber
\end{eqnarray}
and where $N$ will be specified later. Note that $f$ takes values in a
finite-dimensional space, where $Q$ is invertible, by the definition
of the space $\mathcal X$. Recall that $\rho_t(v,x):=D Z^x(t)[v]-
{\mathcal D}_gZ^x(t)$. We have divided the proof into a sequence of lemmas.
\begin{lemma} \label{L2}
We have
%
%
\begin{equation}
\rho_t(v,x) = \zeta^N(v,x)(t) \qquad\forall t\ge0.
\end{equation}
\end{lemma}
\begin{pf}
Adding $f(t)$ to both sides of (\ref{ET211}), we obtain
%
%
\begin{eqnarray}\label{ET215}
&&\frac{\dd\zeta^N(v,x)}{\dd t}(t) +f(t)\nonumber\\
&&\qquad=A\zeta^N(v,x)(t)+
B(Z^x(t),\zeta^N(v,x)(t))\nonumber\\[-8pt]\\[-8pt]
&&\qquad\quad{} +
B(\zeta^N(v,x)(t),Z^x(t)),\nonumber\\
&&\hspace*{-33.25pt}\zeta^N(v,x)(0)=v.\nonumber
\end{eqnarray}
Recall that $D Z^x(t)[v]$ and ${\mathcal D}_gZ^x(t)$ obey equations
(\ref{ET27}) and (\ref{ET28}), respectively. Hence,
$\rho_t:=\rho_t(v,x)$ satisfies
\begin{eqnarray*}
\frac{\dd\rho_t}{\dd t} &=& A\rho_t + B(Z^x(t),\rho_t) +
B(\rho_t,Z^x(t)) - Q^{1/2}g(t),\\
\rho_0 &=& v.
\end{eqnarray*}
Since $f(t)=Q^{1/2}g_t$, we conclude that $\rho_t$ and
$\zeta^N(v,x)(t)$ solve the same linear evolution equation with the
same initial value. Thus, the assertion of the lemma follows.
\end{pf}
\begin{lemma}\label{L3}
For each $N\ge1$, we have $\Pi_{<N}\zeta^N(v,x)(t)=0$ for all $t\ge
2$.
\end{lemma}
\begin{pf} Applying
$\Pi_{<N}$ to both sides of (\ref{ET211}), we obtain
%
%
\begin{eqnarray}\label{ET216}\quad
\frac{\dd}{\dd t}\Pi_{<N}\zeta^N(v,x)(t)&=&-\frac12
\|\Pi_{<N}\zeta^N(v,x)(t)\|_{\mathcal X}^{-1}
\Pi_{<N}\zeta^N(v,x)(t),\nonumber\\[-8pt]\\[-8pt]
\zeta^N(v,x)(0)&=&v.\nonumber
\end{eqnarray}
Multiplying both sides of (\ref{ET216}) by $\Pi_{<N}\zeta^N(v,x)(t)$,
we obtain that $z(t):=\|\Pi_{<N}\zeta^N(v,x)(t)\|_{\mathcal X}^2$ satisfies
%
%
\begin{equation}\label{ode}
\frac{\dd z}{\dd t}(t)=-\frac12 \sqrt{z(t)}.
\end{equation}
Since $\|v\|_{\mathcal X}\le1$, $z(0)\in(0,1]$ and the desired conclusion
holds from elementary properties of the solution of the ordinary
differential equation (\ref{ode}).
\end{pf}
\begin{lemma}\label{L4}
For any $R>0$, the following hold:

\textup{(i)} for any $N$,
%
%
\begin{equation}\label{083102}
\sup_{\|v\|_{\mathcal X}\le1}\sup_{\|x\|_{\mathcal X}\le R}\sup
_{t\in[0,2]} \mathbb{E}
\|\zeta^N(v,x)(t)\|_{\mathcal X}^4 <\infty;
\end{equation}

\textup{(ii)} there exists an $N_0\in\mathbb{N}$ such that for any $N\ge N_0$,
%
%
\begin{equation}\label{083007}
\sup_{\|v\|_{\mathcal X}\le1}\sup_{\|x\|_{\mathcal X}\le R}\int
_0^\infty(
\mathbb{E} \|\zeta^{N}(v,x)(t)\|^4_{\mathcal X} )^{1/2}\,\dd t
<\infty
\end{equation}
and
%
%
\begin{equation}\label{083006}
\sup_{\|v\|_{\mathcal X}\le1}\sup_{\|x\|_{\mathcal X}\le R}\sup
_{t\ge0} \mathbb{E}
\|\zeta^{N}(v,x)(t)\|_{\mathcal X}^4 <\infty.
\end{equation}
\end{lemma}

Since the proof of the lemma is quite lengthy and technical, we
postpone its presentation until the next section. However, we can now
complete the proof of Proposition \ref{m-lm}.

First, we assume that $f$ is given by (\ref{ET213}) with an
arbitrary $N\ge N_0$, where $N_0$ appears in the formulation of
Lemma \ref{L4}. By Lemma \ref{L2}, $\rho_t(v,x)=\zeta^N(v,x)(t)$. Of
course, (\ref{083006}) implies (\ref{ET210}). We show (\ref{ET29}).
As a consequence of Lemma \ref{L3}, we have
$\Pi_{N}\zeta^N(v,x)(t)=0$ for $t\ge2$. The definition of the form
$B(\cdot,\cdot)$ [see (\ref{form-B})] and the fact that the partial
derivatives commute with the projection operator $\Pi_{<N}$ together
imply that
\[
\Pi_{<N} B(Z^x(t), \zeta^N(v,x)(t))= B(Z^x(t),
\Pi_{<N}\zeta^N(v,x)(t)).
\]
As a consequence of Lemma \ref{L3} and convention (\ref{ET212}), we
conclude from (\ref{ET213}) that
\[
f(t)=\Pi_{<N} B(\zeta^N(v,x)(t),Z^x(t))=
B(\zeta^N(v,x)(t),\Pi_{<N} Z^x(t)) \qquad\forall t\ge2.
\]
By (\ref{EA1}), for $t\ge2$, we have
\begin{eqnarray*}
\|f(t)\|_{\mathcal X}&\le& C \|\zeta^N(v,x) (t)\|_{\mathcal X} \|\Pi
_{<N}
Z^x(t)\|_{\mathcal X^{m+1}} \\
&\le& CN \|\zeta^N(v,x) (t)\|_{\mathcal X} \|\Pi_{<N}Z^x(t)\|
_{\mathcal X}
\\
&\le& CN \|\zeta^N(v,x) (t)\|_{\mathcal X} \|Z^x(t)\|_{\mathcal X}.
\end{eqnarray*}
Consequently,
\begin{eqnarray*}
&&\mathbb{E} \int_2^\infty\|Q^{1/2}g_t(v,x)\|_{\mathcal X}^2\,\dd t
\\
&&\qquad=
\mathbb{E} \int_2^\infty\|f(t)\|_{\mathcal X}^2\,\dd t \\
&&\qquad\le C^2N^2 \sup_{t\ge2} (\mathbb{E}
\|Z^x(t)\|_{\mathcal X}^4 )^{1/2}\int_2^\infty( \mathbb
{E}
\|\zeta^N(v,x)(t)\|^4_{\mathcal X} )^{1/2}\,\dd t.
\end{eqnarray*}
Hence, by (\ref{ET23}) and (\ref{083007}), we obtain
\[
\sup_{\|x\|_{\mathcal X}\le R,\|v\|_{\mathcal X}\le1} \mathbb{E}
\int_2^\infty
\|Q^{1/2}g_t(v,x)\|_{\mathcal X}^2\,\dd t <\infty.
\]
Clearly, by Lemma \ref{L4}(i) and (\ref{ET23}),
we have
\[
\sup_{\|x\|_{\mathcal X}\le R,\|v\|_{\mathcal X}\le1} \mathbb{E}
\int_0^2
\|Q^{1/2}g_t(v,x)\|_{\mathcal X}^2\,\dd t <\infty
\]
and the proof of (\ref{ET29}) is completed.

\subsubsection[Proof of Lemma 9]{Proof of Lemma \protect\ref{L4}}

Recall that for any $r$, $A$ is a self-adjoint operator when
considered on the space $\mathcal X^r$, and that
%
%
\begin{equation}
\label{spec-gap} \langle A\psi,\psi\rangle_{\mathcal X^r} \le
-\gamma_*\|\psi\|_{\mathcal X^r}^2,\qquad \psi\in D(A),
\end{equation}
where $\gamma_*>0$ was defined in (\ref{lgamma}). Recall that $V^x$
is the solution to the Ornstein--Uhlenbeck equation (\ref{EOU})
starting from $x$ and that $\mathbf x_x$ is the corresponding solution to
(\ref{Ex1}). The laws of the processes $ (Z^x(t) )_{ t\ge0}$
and $ (V^x(t,\cdot+\mathbf x_x(t)) )_{ t\ge0}$ are the same.
By virtue of this and the fact that $\|V^x(t,\cdot+\mathbf x_x(t))\|
_{\mathcal X}=
\|V^x(t)\|_{\mathcal X}$, we obtain that for each $N\ge1$ and $r\ge0$,
%
%
\begin{equation}
\label{laws} {\mathcal L} ( (\|\Pi_{\ge N} V^x(t)\|
_{\mathcal X
^r} )_{
t\ge0} )= {\mathcal L} ( (\|\Pi_{\ge N} Z^x(t)\|
_{\mathcal X
^r} )_{ t\ge
0} ),
\end{equation}
where, as we recall, ${\mathcal L}$ stands for the law of the respective
process.

In order to show the first part of the lemma, note that from
(\ref{ET211}), upon scalar multiplication (in $\mathcal X$) of both sides
by $\zeta^N(v,x)(t)$ and use of (\ref{spec-gap}),
we have
\begin{eqnarray*}
&&\frac{1}{2} \frac{\dd}{\dd t} \|\zeta^N(v,x)(t)\|_{\mathcal X}^2 \\
&&\qquad\le
-\gamma_*\|\zeta^N(v,x)(t)\|^2_{\mathcal X}\\
&&\qquad\quad{} + |\zeta^N(v,x)(t,0)| \| Z^x(t)\|_{\mathcal V}\|\zeta
^N(v,x)(t)\|
_{\mathcal X} + \frac12 \|\zeta^N(v,x)(t)\|_{\mathcal X}.
\end{eqnarray*}
Here, as we recall, $\mathcal V= \mathcal X^{m+1}$. Taking into account
(\ref{EA1})
and the rough estimate $a/2\le1+a^2$, we obtain
\begin{eqnarray*}
\frac{1}{2} \frac{\dd}{\dd t} \|\zeta^N(v,x)(t)\|_{\mathcal X}^2
&\le&
(1-\gamma_*) \|\zeta^N(v,x)(t)\|^2_{\mathcal X}\\
&&{}+
C\|\zeta^N(v,x)(t)\|^2_{\mathcal X} \| Z^x(t)\|_{\mathcal V} +1.
\end{eqnarray*}
Using Gronwall's inequality and (\ref{laws}), we obtain
\begin{eqnarray*}
\|\zeta^N(v,x)(t)\|^2_{\mathcal X} &\le& (\|v\|^2_{\mathcal
X}+t )
\exp\biggl\{2(1-\gamma_*) t + 2 C \int_0^t \|Z^x(s)\|_{\mathcal V}
\,\dd
s \biggr\}
\\
&\le& (1+t )\exp\biggl\{2(1-\gamma_*) t + C\int_0^t
\|V^x(s)\|_{\mathcal V} \,\dd s \biggr\}\\
&\le&
(1+t )\exp\biggl\{2(1-\gamma_*) t + C\int_0^t
\bigl(\|S(s)x\|_{\mathcal V} + \|M(s)\|_{\mathcal V} \bigr) \,\dd
s \biggr\},
\end{eqnarray*}
where $M=V^0$ is given by (\ref{082903}). By Lemma \ref{L1}(ii),
\[
\sup_{\|x\|_{\mathcal X}\le1}\int_0^\infty\|S(s)x\|_{\mathcal V}\,\dd
s <\infty.
\]
Thus, the proof of the first part of the lemma will be completed as
soon as we can show that
%
%
\begin{equation}\label{Pitu}
\mathbb{E} \exp\biggl\{ C\int_0^2 \|M(s)\|_{\mathcal V} \,\dd s
\biggr\}<\infty.
\end{equation}
By Lemma \ref{L1}, $M$ is a Gaussian element in $C([0,2],\mathcal V)$.
Therefore, (\ref{Pitu}) is a direct consequence of the Fernique theorem
(see, e.g., \cite{DZ1}).

To prove the second part of the lemma, first observe that for
any $N\ge1$,
\begin{eqnarray*}
&&\langle\Pi_{ \ge N} B(Z^x(t),\zeta^N(v,x)(t)),
\zeta^N(v,x)(t)\rangle_{\mathcal X}\\
&&\qquad=\langle B(Z^x(t), \Pi_{ \ge N}\zeta^N(v,x)(t)),
\Pi_{\ge N}\zeta^N(v,x)(t)\rangle_{\mathcal X}\\
&&\qquad=0.
\end{eqnarray*}
Multiplying both sides of (\ref{ET211}) by $\zeta^N(v,x)(t)$ and
remembering that $\Pi_{<N}\zeta^N(v,\break x)(t)=0$ for $t\ge2$, we obtain
that, for those times,
\begin{eqnarray*}
&&\frac{1}{2} \frac{\dd}{\dd t} \|\zeta^N(v,x)(t)\|_{\mathcal X}^2
\\
&&\qquad\le
-\gamma_*\|\zeta^N(v,x)(t)\|^2_{\mathcal X} + |\zeta^N(v,x)(t,0)| \|
\Pi_{\ge N} Z^x(t)\|_{\mathcal V}\|\zeta^N(v,x)(t)\| _{\mathcal X}
\\
&&\hspace*{-4.21pt}\qquad
\stackrel{(\fontsize{8.36}{10}\selectfont{\mbox{\ref{EA1}}})}{\le}-\gamma_*\|\zeta^N(v,x)(t)\|
^2_{\mathcal X} + C
\| \Pi_{\ge N} Z^x(t)\|_{\mathcal V}\|\zeta^N(v,x)(t)\| _{\mathcal
X}^2\\
&&\qquad \le-\gamma_*\|\zeta^N(v,x)(t)\|^2_{\mathcal X}\\
&&\qquad\quad{} + C \bigl(\|
\Pi
_{\ge
N}S(t)x\|_{\mathcal V}^2 + \| \Pi_{\ge N}
Z^0(t)\|_{\mathcal V} \bigr)\|\zeta^N(v,x)(t)\| _{\mathcal X}^2.
\end{eqnarray*}
Define
\[
h(z)= \frac{z^2}{\sqrt{1+\gamma_{*}^{-1}|z|^2}},\qquad z\ge0.
\]
Note that there exists a constant $\tilde C$ such that
\[
Cz\zeta^2 \le\frac{\gamma_{*}}{2}\zeta^2 + \frac{\tilde C}{4}
h(z)\zeta^2,\qquad z\ge0, \zeta\in\mathbb{R}.
\]
Therefore,
\begin{eqnarray*}
&&
C \| \Pi_{\ge N} Z^0(t)\|_{\mathcal V}\|\zeta^N(v,x)(t)\| _{\mathcal
X}^2\\
&&\qquad\le
\frac{\gamma_*}{2}\|\zeta^N(v,x)(t)\| _{\mathcal X}^2 + \frac
{\tilde
C}{4} \|\zeta^N(v,x)(t)\| _{\mathcal X}^2 h(\|\Pi_{\ge N}Z^0(t)\|
_{\mathcal V}).
\end{eqnarray*}
Using Gronwall's inequality, we obtain, for $t\ge2$,
\begin{eqnarray*}
&&\|\zeta^N(v,x)(t)\|^2_{\mathcal X} \\
&&\qquad\le\|\zeta^N(v,x)(2)\|^2_{\mathcal X} \exp\biggl\{-\gamma_* (t-2)
+ L
\|x\|_{\mathcal X}^2\\
&&\qquad\quad\hspace*{93.2pt}{} + \frac{\tilde C}{2} \int_2^t h (\|\Pi_{\ge
N}Z^0(s)\|_{\mathcal V} ) \,\dd s \biggr\},
\end{eqnarray*}
where $ L:=2C \int_2^\infty\|S(t)\|_{L(\mathcal X,\mathcal V)}\,\dd t$.
We have,
therefore, by virtue of the Cauchy--Schwarz inequality,
\begin{eqnarray*}
&&\mathbb{E} \|\zeta^N(v,x)(t)\|_{\mathcal X}^{4}\\
&&\qquad \le\mathbb{E} \|\zeta^N(v,x)(2)\|^4_{\mathcal X}
\mathbb{E} \exp\biggl\{-2\gamma_* (t-2) + 2L \|x\|_{\mathcal
X}^2\\
&&\qquad\quad\hspace*{107.52pt}{} + \tilde
C\int_2^t h (\|\Pi_{\ge N}Z^0(s)\|_{\mathcal V} ) \,\dd
s \biggr\}
\\
&&\hspace*{-5.87pt}\qquad \stackrel{(\fontsize{8.36}{10}\selectfont{\mbox{\ref{laws}}})}{=}
\mathbb{E} \|\zeta^N(v,x)(2)\|
^4_{\mathcal X}
\mathbb{E} \exp\biggl\{-2\gamma_* (t-2) + 2L \|x\|_{\mathcal V}^2
\\
&&\qquad\quad\hspace*{114.4pt}{} + \tilde
C\int_2^t h (\|\Pi_{\ge N}M(s)\|_{\mathcal V} ) \,\dd s
\biggr\},
\end{eqnarray*}
where $M$ is given by (\ref{082903}). Write
\[
\Psi_N(t) :=\exp\biggl\{\tilde C\int_0^t h (\|\Pi_{\ge
N}M(s)\|_{\mathcal V} )\,\dd s \biggr\}.
\]
The proof of part (ii) of the lemma will be completed as soon as we
can show that there exists an $N_0$ such that, for all $N\ge N_0$,
\[
\sup_{t\ge0} e^{-4\gamma_*(t-2)} \mathbb{E} \Psi
_{N}(t)<\infty
\quad\mbox{and}\quad
\int_2^\infty e^{-2\gamma_*t} (\mathbb{E} \Psi
_{N}(t)
)^{1/2}\,\dd t <\infty.
\]
To do this, it is enough to show that
%
%
\begin{equation}
\label{083110} \forall \kappa>0, \exists N(\kappa)\ge1
\qquad\mbox{such that }
\sup_{t\ge0}e^{-\kappa t} \mathbb{E} \Psi
_{N}(t)<\infty\qquad
\forall N\ge N(\kappa).\hspace*{-37pt}
\end{equation}

To do this, note that for any $N_1>N$, $M_{N,N_1}(t):=\Pi_{<N_1}
\Pi_{\ge N}M(t)$ is a strong solution to the equation
\[
\dd M_{N,N_1}(t)= A_N M_{N,N_1}(t) \,\dd t + \Pi_{<M}\Pi_{\ge
N} Q^{1/2} \dd W(t).
\]
Therefore, we can apply the It\^o formula to $M_{N,N_1}(t)$ and
the function
\[
H(x) = \bigl(1+\|(-A_N)^{-1/2}x\|_{\mathcal V}^2 \bigr)^{1/2}.
\]
As a result, we obtain
\begin{eqnarray*}
&& \bigl( 1+ \|(-A_N)^{-1/2} M_{N,N_1}(t)\|_{\mathcal V}^2
\bigr)^{1/2} \\
&&\qquad= 1- \int_0^t \|M_{N,N_1}(s)\|^2 _{\mathcal V} \bigl( 1+
\|(-A_N)^{-1/2}M_{N,N_1}(s)\|_{\mathcal V}^2 \bigr)^{-1/2}\,\dd s \\
&&\qquad\quad{} +\frac{1}{2} \int_0^t \|
H''((-A_N)^{-1/2}M_{N,N_1}(s)) \Pi_{<N_1}\Pi_{\ge N}
Q^{1/2}\|^2_{L_{(\mathrm{HS})}(\mathcal X,\mathcal V)} \,\dd s\\
&&\qquad\quad{} + \int_0^t \bigl( 1+
\|(-A_N)^{-1/2}M_{N,N_1}(s)\|_{\mathcal V}^2 \bigr)^{-1/2}
\\
&&\qquad\quad\hspace*{23.8pt}{}\times
\langle
(-A_N)^{-1}M_{N,N_1}(s) ,
\Pi_{<N_1}\Pi_{\ge N}Q^{1/2}\dd W(s) \rangle_{\mathcal V}.
\end{eqnarray*}
Taking into account the spectral gap property of $A$, we obtain
\begin{eqnarray*}
&&\|M_{N,N_1}(s)\|^2 _{\mathcal V} \bigl( 1+
\|(-A_N)^{-1/2}M_{N,N_1}(s)\|_{\mathcal V}^2 \bigr)^{-1/2}
\\
&&\qquad\ge
\|M_{N,N_1}(s)\|^2 _{\mathcal V} \bigl( 1+\gamma_*^{-1}
\|M_{N,N_1}(s)\|_{\mathcal V}^2 \bigr)^{-1/2}
\\
&&\qquad= h (\|M_{N,N_1}(s)\|_{\mathcal V} ).
\end{eqnarray*}
Therefore,
\[
\tilde C\int_0^t h ( \|M_{N,N_1}(s)\|_{\mathcal V} ) \,\dd s
\le
\tilde C
+ \mathcal M_{N,N_1}(t) + \mathcal R_{N,N_1}(t),
\]
where
\begin{eqnarray*}
\mathcal M_{N,N_1}(t) &=& \tilde C\int_0^t \bigl( 1+
\|(-A_N)^{-1/2}M_{N,N_1}(s)\|_{\mathcal V}^2 \bigr)^{-1/2}\\
&&\hspace*{20.59pt}{}\times\langle
(-A_N)^{-1}M_{N,N_1}(s) ,
\Pi_{<N_1}\Pi_{\ge N}Q^{1/2}\,\dd W(s) \rangle_{\mathcal V}\\
&&{}- \frac{\tilde C^2}{2}\int_0^t \bigl(1+
\|(-A_N)^{-1/2}M_{N,N_1}(s)\|_{\mathcal V}^2\bigr)^{-1}
\\
&&\qquad\hspace*{18.7pt}{}\times\| Q^{1/2}
(-A_N)^{-1}M_{N,N_1}(s)\|^2_{\mathcal V}\,\dd s
\end{eqnarray*}
%
and
\begin{eqnarray*}
\mathcal R_{N,N_1}(t) &=& \frac{\tilde C}{2} \int_0^t \|
H''((-A_N)^{-1/2}M_{N,N_1}(s)) \Pi_{<N_1}\Pi_{\ge N}
Q^{1/2}\|^2_{L_{(\mathrm{HS})}(\mathcal X,\mathcal V)} \,\dd s\\
&&{} + \frac{\tilde C^2}{2}\int_0^t \bigl(1+
\|(-A_N)^{-1/2}M_{N,N_1}(s)\|_{\mathcal V}^2\bigr)^{-1}\\
&&\hspace*{41.1pt}{}\times\|Q^{1/2}
(-A_N)^{-1}M_{N,N_1}(s)\|_{\mathcal V}^2\,\dd s.
\end{eqnarray*}
Since
\begin{eqnarray*}
H''(x)&=& \bigl(1+\|(-A_N)^{-1/2}x\|^2_{\mathcal V} \bigr)^{-1/2}
(-A_N)^{-1}\\
&&{} - \bigl(1+\|(-A_N)^{-1/2}x\|_{\mathcal V}^2 \bigr)^{-3/2}
(-A_N)^{-1}x\\
&&\hspace*{12pt}{}\otimes(-A_N)^{-1}x,
\end{eqnarray*}
there exists a constant $C_1$ such that for all $N$, $N_1$ and $t$,
\begin{eqnarray*}
\mathcal R_{N,N_1}(t) &\le& tC_1 \bigl( \| (-A_N)^{-1}\|_{L(\mathcal
X,\mathcal V)}^4
\|Q^{1/2}\|_{L(\mathcal X,\mathcal V)}^2\\
&&\hspace*{32.3pt}{} +
\|(-A_N)^{-1}Q^{1/2}\|_{L_{(\mathrm{HS})}(\mathcal X,\mathcal V)}^2 \bigr)
\end{eqnarray*}
for all $N_1>N$. Let $\kappa>0$. We can choose sufficiently large
$N_0$ such that for $N\ge N_0$,
\[
C_1 \bigl( \| (-A_N)^{-1}\|_{L(\mathcal X,\mathcal V)}^4
\|Q^{1/2}\|_{L(\mathcal X,\mathcal V)}^2 +
\|(-A_N)^{-1}Q^{1/2}\|_{L_{(\mathrm{HS})}(\mathcal X,\mathcal V)}^2 \bigr) \le
\kappa.
\]
Since $ (\exp\{\mathcal M_{N_0,N_1}(t)\} )$ is a
martingale, we
have shown, therefore, that for $N\ge N_0$,
\[
\mathbb{E}\int_0^t \exp\{h (\|M_{N,N_1}(s)\|_{\mathcal
V} ) \,\dd s
\}\le\exp\{ \tilde C + \kappa t \}.
\]
Letting $N_1\to\infty$, we obtain (\ref{083110}).

\section[Proof of Theorem 4]{Proof of Theorem \protect\ref{TTracer}}\label{sec5}
With no loss of generality, we will assume that the initial position
of the tracer $\mathbf x_0=0$.
By definition,
\[
\mathbf x(t)=\int_0^tV(s,\mathbf x(s))\,\dd s =\int_0^t \mathcal
Z(s,0)\,\dd s,
\]
where $\mathcal Z(t,x)= V(t,\mathbf x(t)+x)$ is the observation
process. Recall
that $\mathcal Z(t)$ is a stationary solution to (\ref{Zet}). Obviously,
uniqueness and the law of a stationary solution do not depend on the
particular choice of the Wiener process. Therefore,
\[
\mathcal L \biggl(\frac{\mathbf x(t)}{t} \biggr) = \mathcal L
\biggl(\frac1 t \int_0^t
\tilde Z(s,0)\,\dd s \biggr)
\]
and
\[
\mathcal L \biggl(\frac
{\dd
\mathbf x}{\dd t}(t) \biggr) = \mathcal L (\tilde Z(t,0) ),
\]
where, as before, $\mathcal L(X)$ stands for the law of a random
element $X$
and $\tilde Z$ is, by Theorem \ref{Pttheorem}, a unique (in law)
stationary solution of the equation
\[
\dd\tilde Z(t) = [A \tilde Z(t) + B(\tilde Z(t),\tilde
Z(t)) ]\,\dd t + {Q}^{1/2} \,\dd
W(t).
\]
Let $F\dvtx\mathcal X\to\mathbb{R}$ be given by $F(x)= x(0)$. The
proof of
the first part of the theorem will be completed as soon as we can
show that the limit (in probability)
\[
\mathop{{\mathbb{P}} \mbox{-}\mathrm{lim}}_{t\uparrow\infty} \frac
1 t
\int_0^t \tilde Z(s,0)\,\dd s
\]
exists and is equal to $\int_{\mathcal X} F(x)\mu_*(\dd x)$, where
$\mu_*$
is the unique invariant measure for the Markov family $Z$ defined by
(\ref{EW2}). Since the semigroup $(P_t)_{t\ge0}$ satisfies the
e-property and
is weak-$^*$ mean ergodic, part (2) of Theorem \ref{Theorem}
implies that for any bounded Lipschitz continuous function $\psi$,
\[
\mathop{{\mathbb{P}} \mbox{-}\mathrm{lim}} _{t\uparrow
\infty}\frac1 t \int_0^t \psi( \tilde Z(s))\,\dd s
=\int_{\mathcal X}\psi(x)\mu_*(\dd x).
\]
Since $\mathcal X$ is embedded in the space of bounded continuous
functions, $F$ is Lipschitz. The theorem
then follows by an easy truncation argument.

\section*{Acknowledgments}
The authors wish to express their
gratitude to an anonymous referee for thorough reading of the
manuscript and valuable remarks. We also would like to express our
thanks to Z. Brze\'{z}niak for many enlightening discussions on the
subject of the article.


%
\printaddresses

\end{document}